\documentclass[a4paper,8pt,fleqn]{article}

\usepackage[a4paper,top=0.7in,bottom=0.8in,left=.8in,right=.8in]{geometry}

\usepackage{color}
\usepackage{graphicx}

\definecolor{egyptianblue}{rgb}{0.06,0.2,0.65}
\usepackage{tikz}
\usepackage[tickmarkheight=.5em,textwidth=7em,textsize=tiny,backgroundcolor=gray!25,linecolor=gray!70,bordercolor=gray!70, textcolor=black]{todonotes}

\usepackage{caption}
\usepackage{booktabs}

\usepackage[colorlinks=true,allcolors=egyptianblue,bookmarks=false,hypertexnames=true]{hyperref}

\usepackage{newtxtext}
\usepackage[libertine]{newtxmath}
\usepackage[cal=boondoxo,bb=boondox,frak=boondox]{mathalfa}
\usepackage[scaled=.95,type1]{cabin}
\usepackage{fontawesome} 

\usepackage{adforn}

\usepackage{amsthm}

\usepackage{fancyhdr}
\usepackage{titlesec}
\usepackage{authblk}
\titleformat{\section}[hang]{\large\normalfont\scshape\filcenter}{\thesection.}{1em}{}
\titleformat{\subsection}[hang]{\normalfont\itshape\filcenter}{\thesubsection.}{1em}{}
\titleformat{\paragraph}[runin]{\itshape}{\theparagraph.}{1em}{}

\newtheorem{theorem}{Theorem}[section]
\newtheorem{lemma}[theorem]{Lemma}
\newtheorem{proposition}[theorem]{Proposition}
\newtheorem{corollary}[theorem]{Corollary}

\newtheorem{remark}[theorem]{Remark}

\usepackage{doi}
\usepackage{natbib}
\begin{document}
 \title{\bfseries Noise-resilient penalty operators based on statistical \\ differentiation schemes}
   \author[1,2]{Marc Vidal \thanks{\textbf{Funding}: The authors acknowledge the support from the FWO project G013024N.}
   }
   \author[1]{Yves Rosseel}
   
   \affil[1]{\scshape Department of Data Analysis, Henri Dunantlaan 1, 9000 Ghent, Belgium}
   \affil[2]{\scshape Max Planck Institute for Human Cognitive and Brain Sciences, Leipzig, Germany}
   \affil[ ]{\small\ttfamily
              \href{mailto:marc.vidalbadia@ugent.be}{marc.vidalbadia@ugent.be};
              \href{mailto:yves.rosseel@ugent.be}{yves.rosseel@ugent.be}}
   \date{}

\maketitle
\thispagestyle{empty}
\setcounter{page}{0}

\begin{abstract}
\noindent Classical approaches to penalized smoothing often rely on basis or kernel
expansions, which constrain the estimator to a fixed span and may be
restrictive for discretely observed data. We instead regularize a single
noisy trajectory directly on its observation grid, using difference
operators that remain genuine finite-difference approximations to
derivatives while being statistically normalized and mutually decorrelated
under a reference noise law. We extend this white-noise construction to a
parametric family of covariance-adapted reference geometries, with the
parameter estimated by generalized method of moments from an independent
calibration sample. A first-order plug-in expansion shows that, within this
moment family, efficient calibration minimizes the leading-order
discrepancy between the plug-in and oracle covariance-adapted smoothers.
Numerical experiments confirm the predicted rate at which this plug-in
discrepancy vanishes and compare the resulting reconstruction against
conventional discrete, basis, and kernel smoothers.
\end{abstract}

\vspace{10pt}
{ \small \noindent  {\bf MSC 2020 subject classifications:} Primary 62G08; secondary 62G20. }



\section{Introduction}
\label{sec:1}

Let $X=(X_1,\ldots,X_d)$ be a random vector in $\mathbb R^d$, with
coordinates indexed by $t\in\mathcal I_d=\{1,\ldots,d\}$, which may represent
time, spatial location, or another ordered structure. We interpret an observed realization
$X_i=(X_{i1},\ldots,X_{id})$ as a noisy discretization of an underlying curve
and impose regularity directly on the observed grid, without invoking
continuous function spaces or classical smoothness assumptions. We define
the penalized estimate $\widehat X_i$ as the solution to
\begin{equation}
\label{pen:gen}
\widehat X_i(R;\boldsymbol\alpha)
=
\arg\min_{f\in\mathbb R^d}
\frac12
\left\{
\|X_i-f\|^2+
\sum_{r=1}^R\alpha_r\|D^{(r)}f\|^2
\right\},
\end{equation}
where $D^{(r)}$ is a discrete difference operator of order $r$ and
$\alpha_r>0$ controls the strength of penalization at that order. Each term
$\|D^{(r)}f\|^2$ measures roughness at a different discrete derivative order,
so that \eqref{pen:gen} combines several scales of local variation. The
statistical problem is then not only how strongly these components should be
penalized, but also how roughness measured at different orders should be
calibrated relative to one another.

\begin{figure}[t]
\centering
\includegraphics[width=\textwidth]{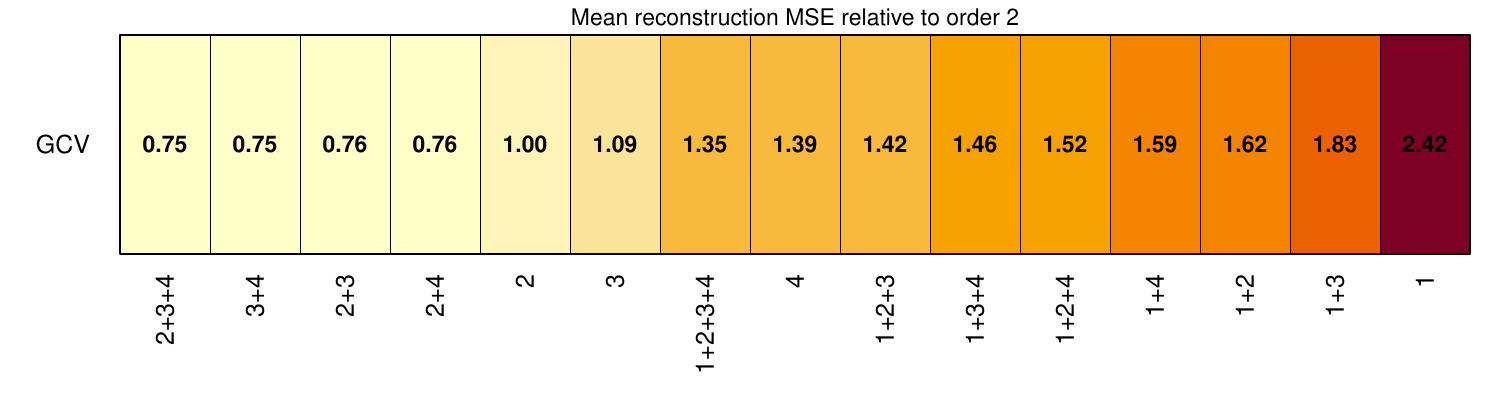}
\caption{Relative reconstruction error for different combinations of
difference orders. A sinusoidal signal is observed at $d=100$ equally spaced
grid points with independent Gaussian noise of standard deviation $0.2$, over
500 Monte Carlo replications. For each nonempty subset
$A\subseteq\{1,2,3,4\}$, the penalty
$\sum_{r\in A}(\bar D^{(r)})^\top\bar D^{(r)}$ is formed using the
uncorrelated difference operators introduced in Section~\ref{sec:2}, with a
common regularization parameter selected by GCV. Values are mean
reconstruction errors relative to the second-order penalty, which is
normalized to one; hence values below one indicate improvement. The results
show that combining selected derivative orders can substantially improve
reconstruction, while including additional orders indiscriminately need not
do so.}
\label{fig:intro-orders}
\end{figure}

Discrete difference-penalized smoothing along these lines has a separate,
older history: minimizing a fidelity term plus a penalty on finite
differences of the raw observations, with no
continuous-function representation. \citet{Whittaker:1923} introduced this
for a fixed third-order penalty in the graduation of mortality tables,
\citet{Hodrick:Prescott:1997} rediscovered the second-order case for
macroeconomic trend extraction, and \citet{Eilers:2003} generalized the
penalty to arbitrary order as the ``perfect smoother''; see also
\citet{Tibshirani:2022}.
The two traditions later reconverged in the P-spline literature, which places a
difference penalty on spline coefficients rather than on the raw grid
\citep{Eilers:Marx:Durban:2015}.
This discrete starting point contrasts with the continuous nonparametric
asymptotic regime in which classical optimality results for penalized
smoothing are usually derived
\citep{Mueller:1987,Hall:Opsomer:2005,Claeskens:Krivobokova:Opsomer:2009}.
Staying on the fixed grid matters because smoothing need not improve
performance: in a related discrete
functional setting, \citet{Belhakem:Picard:Rivoirard:Roche:2025} show that
projection without additional smoothing can already attain the minimax rate,
whereas unsuitable regularization may deteriorate it.
Our motivation is therefore not that \eqref{pen:gen} universally improves
the convergence rate, but that, when discrete regularization is desired,
roughness components of different orders should be measured on a common
statistical scale; see \autoref{fig:intro-orders} for a simple numerical
illustration.

We obtain such a scale from the uncorrelated difference operators introduced
by \citet{Mizuta:2006,Mizuta:2023}. Under the white-noise reference law
$Q_0=\mathcal N(0,I_d)$, finite linear combinations of translates of
conventional differences can be chosen so that distinct orders are
uncorrelated. We show that these operators retain their finite-difference
interpretation: an order-$r$ construction annihilates polynomials of degree
below $r$ and, after the appropriate derivative normalization, provides an
$r$th-order finite-difference approximation. Unit-variance normalization and
cross-order decorrelation are separate properties. The former gives a common
reference scale and, by Parseval's identity, unit average spectral energy at
every order, whereas conventional $r$th differences have energy
$\binom{2r}{r}\sim4^r/\sqrt{\pi r}$. The latter removes the additional
non-negative cross-order term in the variance of multi-order roughness that
remains even after conventional differences are rescaled to unit variance
(see Proposition~\ref{prop:variance-optimality}).

To accommodate dependent reference noise, we subsequently
extend the construction to a parametric family of reference laws
$Q_\vartheta=\mathcal N(0,\Sigma_\vartheta)$,
$\vartheta\in\Theta$, where $\Theta\subset\mathbb R^k$ is an open parameter set.
Here $\Sigma_\vartheta$ defines the covariance geometry used to calibrate the
penalty operators, without imposing a parametric model on the observed curves
$X_i$. Replacing Euclidean orthogonality by the $\Sigma_\vartheta$ inner
product produces covariance-adapted difference operators while leaving
derivative consistency unchanged. Whenever the corresponding local rank
condition holds, these operators are unique up to sign and depend smoothly
on $\vartheta$. This condition can become restrictive at higher orders,
where orthogonality alone need not determine a unique operator. When $\vartheta$ is unknown, the reference covariance geometry must itself
be estimated from an independent calibration sample
$\varepsilon_1,\ldots,\varepsilon_n\overset{\mathrm{iid}}{\sim}Q_{\vartheta_0}$,
drawn from the reference law itself rather than from the curves being
regularized, before the covariance-adapted operators can be implemented.
The standardized multi-order roughness coordinates of the $\varepsilon_i$,
whose expectations vanish at the correctly calibrated geometry, provide
estimating equations for $\vartheta$. Under Gaussianity, their
covariance reduces to $2I_R$ after decorrelation, so that, when the number of
moment conditions exceeds the dimension of $\vartheta$, 
equal weighting is efficient within this family of estimating equations.
For a root-$n$-consistent estimator
$\widehat\vartheta_n$ of $\vartheta_0$, standard generalized method of
moments (GMM) theory provides its asymptotic normal distribution under the usual
identification and regularity conditions. Separately, using the framework of
\citet{Schick:2001}, we obtain an asymptotic linearization of the
parameter-dependent empirical estimating equations evaluated at
$\widehat\vartheta_n$, which is needed for the roughness-based diagnostics
computed at the estimated geometry. Smooth dependence of the
covariance-adapted difference operators on $\vartheta$, together with
finite-dimensional matrix calculus, then propagates this estimation
uncertainty through the penalized reconstruction of any observed curve
$X_i$ independent of the calibration sample. In particular, replacing
the oracle geometry $\vartheta_0$ by $\widehat\vartheta_n$ produces a
first-order perturbation of order $n^{-1/2}$ in the regularized curve
$\widehat X_i=S_{\widehat\vartheta_n}X_i$ and, consequently, an order
$n^{-1}$ squared plug-in discrepancy. Moreover, efficient GMM
estimation within this moment family minimizes the leading-order
plug-in discrepancy relative to the oracle covariance-adapted reconstruction.
Thus, the asymptotic theory quantifies the cost of learning the reference
covariance geometry, from a sample of the reference law itself, before using
it for discrete regularization.

The paper is organized as follows.
Sections~\ref{sec:2} and~\ref{sec:3} introduce the discrete difference operators
and the associated smoothing operators.
Section~\ref{sec:4} extends this white-noise construction to a parametric
family of covariance-adapted reference geometries.
Section~\ref{sec:5} estimates the reference parameter from an independent
calibration sample and propagates the resulting uncertainty into the
smoother.
Section~\ref{sec:6} discusses curve-wise selection of the regularization
strength by generalized cross-validation, and Section~\ref{sec:7} presents
numerical experiments. Proofs of the formal statements are collected in
Appendix~\ref{appn}, while an additional reconstruction experiment is
reported in Appendix~\ref{appB}.

\section{Difference operators with decorrelation properties}
\label{sec:2}

A discrete difference operator of order $r$ is a linear map
$D^{(r)}:\mathbb{R}^d\to\mathbb{R}^{d_r}$, defined by convolution with a finite
stencil of signed weights $w_{-L_r}^{(r)},\ldots,w_{L_r}^{(r)}$. The operator
acts on $x\in\mathbb{R}^d$ as
\[
(D^{(r)}x)(t)=\sum_{\ell=-L_r}^{L_r}w_\ell^{(r)}x_{t+\ell},
\qquad
t\in\mathcal I^{(r)}=\{L_r+1,\ldots,d-L_r\},
\]
where the stencil is fully supported, so that $d_r=d-2L_r$. Classical
finite-difference operators use stencils determined by binomial coefficients
and their shifted or centred versions.

Following \cite{Mizuta:2023}, we instead consider difference operators
constructed to satisfy a statistical decorrelation property. Let
$Q_0\in\mathscr M_1(\mathbb R^d)$ denote the white-noise reference law under
which the entries of a reference vector $\varepsilon=(\varepsilon_1,\ldots,\varepsilon_d)$
are independent standard normal random variables; $\varepsilon$ is used here
purely to define the decorrelation geometry and is distinct from the observed
curves $X_i$ introduced in Section~\ref{sec:1}. Hence $E_{Q_0}(\varepsilon)=0$ and
$\operatorname{cov}_{Q_0}(\varepsilon)=I_d$. For any two finite stencils $u$ and $v$,
applied at the same grid location $t$,
\[
\operatorname{cov}_{Q_0}
\left\{\sum_\ell u_\ell \varepsilon_{t+\ell},
       \sum_\ell v_\ell \varepsilon_{t+\ell}\right\}
=\sum_\ell u_\ell v_\ell,
\]
where stencils with different support widths are understood to be zero-padded
to a common support. Thus decorrelation under $Q_0$ corresponds exactly to
Euclidean orthogonality of the associated stencil vectors.

 \cite{Mizuta:2023} constructs the $r$th uncorrelated difference as a
finite linear combination of translates of the conventional $r$th-order
difference. To make this structure explicit, let $T$ denote the unit forward
shift, $(Tx)_t=x_{t+1}$, and define
\[
\Delta^{(r)}=(1-T^{-1})^r,
\qquad
(\Delta^{(r)}x)_t
=\sum_{k=0}^r(-1)^k\binom{r}{k}x_{t-k}.
\]
The $r$th uncorrelated difference operator can then be written as
$\bar D^{(r)}=W_r(T)\Delta^{(r)}$, where
$W_r(T)=\sum_{j\in J_r}a_j^{(r)}T^j$ for a finite index set
$J_r\subset\mathbb Z$. Equivalently,
\[
(\bar D^{(r)}x)_t
=\sum_{j\in J_r}a_j^{(r)}(\Delta^{(r)}x)_{t+j}.
\]
The coefficients $a_j^{(r)}$ are chosen according to the recursive
construction of \citet{Mizuta:2023}, which enforces mutual decorrelation at
a common grid location under $Q_0$, with minimal support subject to this
constraint and unit-variance normalization. We retain this construction but make explicit below a property not
established in \cite{Mizuta:2023}: statistical decorrelation within this class
of operators preserves their interpretation as finite-difference
approximations of the corresponding derivative order.

The first uncorrelated difference operators given in \cite{Mizuta:2023} are
displayed below. Each row contains the unnormalized stencil coefficients,
centred at $t$.

\medskip

\[
\begin{array}{c@{\quad}c}
&
\begin{tikzpicture}[scale=1,
  stencil/.style={minimum size=1.5em, inner sep=0pt, text height=1.5ex, text depth=.25ex}
]

\foreach \r/\row/\L in {
  0/{0,0,0,0,0,1,0,0,0,0,0}/0,
  1/{0,0,0,0,1,0,-1,0,0,0,0}/1,
  2/{0,0,0,1,-1,0,-1,1,0,0,0}/2,
  3/{0,0,2,-3,0,0,0,3,-2,0,0}/3,
  4/{7,-16,9,0,0,0,0,0,9,-16,7}/5
}{
  \pgfmathtruncatemacro{\left}{5 - \L}
  \pgfmathtruncatemacro{\right}{5 + \L}

  \fill[gray!10,rounded corners]
    (\left - 0.5, -\r - 0.5) rectangle (\right + 0.5, -\r + 0.5);

  \foreach \val [count=\i] in \row {
    \pgfmathtruncatemacro{\x}{\i - 1}
    \pgfmathtruncatemacro{\y}{-\r}
    \pgfmathsetmacro{\valnum}{\val}
    \ifdim\valnum pt=0pt
      \def\col{blue!10}
    \else
      \ifdim\valnum pt>0pt
        \def\col{green!30}
      \else
        \def\col{red!30}
      \fi
    \fi
    \node[draw, fill=\col, font=\footnotesize, stencil] at (\x, \y) {\val};
  }
}

\foreach \label [count=\i from 0] in {
  $\bar{D}^{(0)}$,
  $\bar{D}^{(1)}$,
  $\bar{D}^{(2)}$,
  $\bar{D}^{(3)}$,
  $\bar{D}^{(4)}$
}{
  \node[anchor=east] at (-0.8,-\i) {\label};
}

\foreach \label [count=\i from 0] in {
  $t-5$,$t-4$,$t-3$,$t-2$,$t-1$,
  $t$,
  $t+1$,$t+2$,$t+3$,$t+4$,$t+5$
}{
  \node[anchor=south] at (\i,0.7) {\label};
}

\end{tikzpicture}
\end{array}
\]

\medskip

The construction above is fundamentally different from choosing an arbitrary
stencil and subsequently imposing orthogonality constraints. The statistical
orthogonalization is performed within the class generated by translates of
the conventional $r$th-order difference. Consequently, the algebraic
cancellation properties of an $r$th difference are retained. The following
result makes this precise.

\begin{proposition}
\label{prop:1}
Let $h>0$ denote the grid spacing, so that grid location $i$
corresponds to the real coordinate $t_i=t_0+ih$, and suppose that
$\bar D^{(r)}=W_r(T)\Delta^{(r)}$, where
$W_r(T)=\sum_{j\in J_r}a_j^{(r)}T^j$ and $J_r\subset\mathbb Z$ is finite.
Write the resulting stencil as
$(\bar D^{(r)}x)_i=\sum_\ell w_\ell^{(r)}x_{i+\ell}$ and define its discrete
moments by $M_q(w^{(r)})=\sum_\ell \ell^q w_\ell^{(r)}$. Then
\[
M_q(w^{(r)})=0,\quad q=0,\ldots,r-1,
\qquad\text{and}\qquad
M_r(w^{(r)})=r!\sum_{j\in J_r}a_j^{(r)}.
\]
Hence the resulting stencil has exact difference order $r$ whenever
$\sum_{j\in J_r}a_j^{(r)}\neq0$. Moreover, if $f\in C^{r+1}$ in a
neighbourhood of $t_i$ and $x_i=f(t_i)$, then
\[
(\bar D^{(r)}x)_i
=
h^r\left(\sum_{j\in J_r}a_j^{(r)}\right)f^{(r)}(t_i)
+O(h^{r+1}).
\]
If, in addition, $w_{-\ell}^{(r)}=(-1)^r w_\ell^{(r)}$, then
$M_{r+1}(w^{(r)})=0$ and, for $f\in C^{r+2}$, the remainder improves to
$O(h^{r+2})$.
\end{proposition}

\begin{remark}
\label{rem:derivative-normalization}
The normalization used in the construction of $\bar D^{(r)}$ is statistical:
its stencil is scaled to have unit variance under $Q_0$. This should be
distinguished from the scaling required to approximate the $r$th derivative.
Indeed, Proposition~\ref{prop:1} shows that, for
$c_r=\sum_{j\in J_r}a_j^{(r)}\neq0$,
\[
(\bar D^{(r)}x)_i
=
c_rh^r f^{(r)}(t_i)+O(h^{r+1}),
\qquad x_i=f(t_i).
\]
Thus the unscaled difference operator approximates
$c_rh^r f^{(r)}(t_i)$ rather than $f^{(r)}(t_i)$ itself. Dividing by the
leading factor $c_rh^r$ gives
\[
(c_rh^r)^{-1}(\bar D^{(r)}x)_i
=
f^{(r)}(t_i)+O(h).
\]
Under the parity condition of Proposition~\ref{prop:1}, and for
$f\in C^{r+2}$, the unscaled remainder is $O(h^{r+2})$, so the corresponding
normalized derivative approximation has error $O(h^2)$.
The division by $c_rh^r$ is not an additional normalization imposed in the
construction of $\bar D^{(r)}$. The operators used throughout this paper
retain their unit-variance normalization under $Q_0$; the preceding
rescaling only makes explicit their interpretation as finite-difference
approximations of the corresponding derivatives.
\end{remark}

For the stencils reported by \cite{Mizuta:2023}, these moments can
be checked directly from their coefficients. In particular,
\[
\begin{aligned}
M_1(w^{(1)})&=2(-1)=-2,\\
M_2(w^{(2)})&=2\{1(-1)+4(1)\}=6,\\
M_3(w^{(3)})&=2\{8(3)+27(-2)\}=-60,\\
M_4(w^{(4)})&=2\{81(9)+256(-16)+625(7)\}=2016.
\end{aligned}
\]
Thus the first nonzero discrete moments are
$M_1=-2$, $M_2=6$, $M_3=-60$, and $M_4=2016$, respectively, while all moments
below the corresponding order vanish. Their symmetric or antisymmetric
structure also gives $M_{r+1}=0$ in each case. Hence the concrete operators
proposed by Mizuta satisfy the conditions of Proposition~\ref{prop:1}.
\begin{corollary}
\label{col:1}
Let $\varepsilon\sim Q_0$, where the entries of $\varepsilon$ are independent with zero mean and
unit variance, and let $\bar D^{(0)},\ldots,\bar D^{(R)}$ be uncorrelated
difference operators constructed as above and normalized to unit variance
under $Q_0$. Then, at every grid location $t$ at which the corresponding
stencils are fully supported,
\[
E_{Q_0}\{(\bar D^{(r)}\varepsilon)_t\}=0,\qquad
\operatorname{var}_{Q_0}\{(\bar D^{(r)}\varepsilon)_t\}=1,
\qquad
\operatorname{cov}_{Q_0}
\{(\bar D^{(r)}\varepsilon)_t,(\bar D^{(s)}\varepsilon)_t\}=0
\quad(r\neq s).
\]
If $Q_0=\mathcal N(0,I_d)$, the corresponding transformed variables are
independent across difference orders at each fixed $t$.
\end{corollary}

The decorrelation property in Corollary~\ref{col:1} concerns different
difference orders evaluated at the same grid location. It does not imply
independence across grid locations. Indeed, overlapping stencils generally
induce dependence between $(\bar D^{(r)}\varepsilon)_t$ and
$(\bar D^{(s)}\varepsilon)_{t'}$ when $t\neq t'$, including within a fixed difference
order. The construction therefore removes the artificial cross-order
correlation generated by the difference operators under the reference
white-noise model, rather than attempting to decorrelate the transformed
process across the grid.

The unit-variance normalization established in Corollary~\ref{col:1} also
has a frequency-domain interpretation
that is independent of cross-order decorrelation. Let $w^{(r)}$ denote the
zero-padded stencil of $\bar D^{(r)}$ and write
$\widehat w^{(r)}(\omega_j)=
\sum_\ell w_\ell^{(r)}e^{-i\omega_j\ell}$, where
$\omega_j=2\pi j/d$.

\begin{proposition}
\label{prop:spectral-normalization}
For every normalized difference order $r=1,\ldots,R$,
\[
\frac{1}{d}\sum_{j=0}^{d-1}
\left|\widehat w^{(r)}(\omega_j)\right|^2=1.
\]
Hence all normalized difference orders have the same average spectral energy
under the reference white-noise geometry.
\end{proposition}

This normalization should be distinguished from decorrelation. Indeed, the
conventional $r$th-order difference has stencil coefficients
$(-1)^k\binom{r}{k}$ and therefore squared Euclidean norm
\[
\sum_{k=0}^r\binom{r}{k}^2
=
\binom{2r}{r}
\sim\frac{4^r}{\sqrt{\pi r}}.
\]
Thus conventional differences are naturally expressed on strongly
order-dependent scales. Rescaling each conventional stencil to unit norm
removes this scale effect and gives the same conclusion as
Proposition~\ref{prop:spectral-normalization}, but it does not remove
cross-order dependence.

\section{Discrete smoothing operators}
\label{sec:3}

\subsection{Discrete penalization with finite differences}
We begin by motivating the discrete penalty construction from the usual
continuous roughness-penalization perspective. A common approach is to replace
the $L^2$ inner product by a Sobolev-type inner product of the form
$
\langle f, g \rangle_\alpha = \langle f, g \rangle_{L^2} + \alpha \langle f^{(r)}, g^{(r)} \rangle_{L^2},
$
which formally induces the smoothing operator
$
\mathcal{S}^{(r)}_\alpha = \left(\mathbb{I} + \alpha\, \mathcal{D}^{(r)*} \mathcal{D}^{(r)}\right)^{-1},
$
where $\mathcal{D}^{(r)}$ denotes the $r$th derivative operator, $\mathcal{D}^{(r)*}$ its adjoint, and $\mathbb{I}$ represents the identity in the space.  This operator serves as the continuous prototype for our discrete penalty framework. 

In the discrete setting, let
$D^{(r)}\in\mathbb R^{d_r\times d}$ denote a finite-difference operator of
order $r$. For an observed trajectory $X_i\in\mathbb R^d$, the corresponding
roughness penalty is
\[
X_i^\top\mathfrak P^{(r)}X_i
=
\|D^{(r)}X_i\|^2,
\qquad
\mathfrak P^{(r)}
=
D^{(r)\top}D^{(r)}.
\]
The associated smoothing operator is
\begin{equation}
\label{op:smoothing}
S_\alpha^{(r)}
=
\left(
I_d+\alpha\mathfrak P^{(r)}
\right)^{-1},
\end{equation}
and the penalized estimator
$\hat X_i=S_\alpha^{(r)}X_i$ is the solution of
\[
\hat X_i
=
\arg\min_{f\in\mathbb R^d}
\frac12
\left\{
\|X_i-f\|^2
+
\alpha\|D^{(r)}f\|^2
\right\}.
\]

At this level, $D^{(r)}$ may be any finite-difference operator.
The construction of Section~\ref{sec:2}, however, provides an alternative
penalty
\[
\mathfrak P_{Q_0}^{(r)}
=
\bar D^{(r)\top}\bar D^{(r)},
\]
where $\bar D^{(r)}$ preserves the $r$th finite-difference order while being
normalized to unit variance under the reference white-noise law $Q_0$.
Accordingly,
\[
S_{\alpha,Q_0}^{(r)}
=
\left(
I_d+\alpha\mathfrak P_{Q_0}^{(r)}
\right)^{-1}
\]
penalizes the same derivative order as its conventional counterpart, but using
a statistically normalized difference operator.

For a fixed order $r$, we also consider the interpolating family
\begin{equation}
\label{eq:blendpen}
\mathfrak P_\eta^{(r)}
=
(1-\eta)\mathfrak P^{(r)}
+
\eta\mathfrak P_{Q_0}^{(r)},
\qquad
\eta\in[0,1],
\end{equation}
with corresponding smoother
\[
S_{\alpha,\eta}^{(r)}
=
\left(
I_d+\alpha\mathfrak P_\eta^{(r)}
\right)^{-1}.
\]
Thus $\eta=0$ gives conventional finite-difference penalization and
$\eta=1$ gives the statistically decorrelated construction.

\subsection{Statistical structure of multi-order roughness}

The statistical distinction between conventional and uncorrelated difference
operators becomes particularly transparent when several difference orders are
considered simultaneously. Let $R$ be fixed and consider an interior grid
location $t$ at which all stencils
$\bar D^{(1)},\ldots,\bar D^{(R)}$ are fully supported. Define, for the
reference vector $\varepsilon\sim Q_0$ of Section~\ref{sec:2},
$
Z(t)
=
\bigl(
(\bar D^{(1)}\varepsilon)_t,\ldots,(\bar D^{(R)}\varepsilon)_t
\bigr)^\top.
$
By Corollary~\ref{col:1}, under the reference law
$Q_0=\mathcal N(0,I_d)$,
$
Z(t)\sim\mathcal N(0,I_R).
$
Thus the construction provides two distinct properties: the coordinates
associated with different difference orders are mutually orthogonal under the
reference law, and each coordinate is normalized to have unit variance. The
first property concerns cross-order dependence, whereas the second places the
different orders on a common statistical scale.

For non-negative order-specific weights
$\boldsymbol\alpha=(\alpha_1,\ldots,\alpha_R)^\top$, define the local
multi-order roughness statistic
\[
\mathcal R_{\boldsymbol\alpha}(\varepsilon;t)
=
\sum_{r=1}^R
\alpha_r(\bar D^{(r)}\varepsilon)_t^2.
\]
Its reference distribution is available exactly.

\begin{proposition}
\label{prop:multiorder}
Let $\varepsilon\sim Q_0=\mathcal N(0,I_d)$ and let
$\bar D^{(1)},\ldots,\bar D^{(R)}$ be uncorrelated difference operators
normalized as in Corollary~\ref{col:1}. Then, at every grid location $t$ at
which all stencils are fully supported,
\[
\mathcal R_{\boldsymbol\alpha}(\varepsilon;t)
\overset{d}{=}
\sum_{r=1}^R\alpha_r\chi_{1,r}^2,
\]
where $\chi_{1,1}^2,\ldots,\chi_{1,R}^2$ are independent chi-squared random
variables with one degree of freedom. Consequently,
\[
E_{Q_0}\{\mathcal R_{\boldsymbol\alpha}(\varepsilon;t)\}
=\sum_{r=1}^R\alpha_r,
\qquad
\operatorname{var}_{Q_0}\{\mathcal R_{\boldsymbol\alpha}(\varepsilon;t)\}
=2\sum_{r=1}^R\alpha_r^2.
\]
\end{proposition}

To isolate the effect of decorrelation specifically, consider any family of
difference stencils normalized to unit variance under $Q_0$, not necessarily
mutually decorrelated. At a fixed location, let
$Z_G(t)\sim\mathcal N(0,G)$ denote the corresponding vector of difference
coordinates, where $G$ is their correlation matrix and hence $G_{rr}=1$.
For $A=\operatorname{diag}(\alpha_1,\ldots,\alpha_R)$, define
$\mathcal R_{\boldsymbol\alpha,G}(\varepsilon;t)=Z_G(t)^\top A Z_G(t)$.
The usual Gaussian quadratic-form identities give
\[
E_{Q_0}\{\mathcal R_{\boldsymbol\alpha,G}(\varepsilon;t)\}
=\operatorname{tr}(AG),
\qquad
\operatorname{var}_{Q_0}\{\mathcal R_{\boldsymbol\alpha,G}(\varepsilon;t)\}
=2\operatorname{tr}\{(AG)^2\}.
\]
Since $G_{rr}=1$, the expectation is $\sum_r\alpha_r$ for every such
normalized family. The variance, however, retains the effect of cross-order
dependence.

\begin{proposition}
\label{prop:variance-optimality}
Let $G$ be the correlation matrix under $Q_0$ of any family of
unit-variance difference coordinates, and let
$A=\operatorname{diag}(\alpha_1,\ldots,\alpha_R)$ with $\alpha_r>0$.
Then
\[
\operatorname{var}_{Q_0}
\{\mathcal R_{\boldsymbol\alpha,G}(\varepsilon;t)\}
=
2\sum_{r=1}^R\alpha_r^2
+
2\sum_{r\neq s}
\alpha_r\alpha_sG_{rs}^2.
\]
Consequently,
$\operatorname{var}_{Q_0}\{\mathcal R_{\boldsymbol\alpha,G}(\varepsilon;t)\}
\ge 2\sum_r\alpha_r^2$, with equality if and only if $G=I_R$.
Thus, among unit-variance multi-order representations with fixed positive
weights $\boldsymbol\alpha$, mutual decorrelation uniquely minimizes the
reference variance of the aggregate roughness statistic.
\end{proposition}

Proposition~\ref{prop:variance-optimality} separates the role of
decorrelation from that of normalization. Normalizing conventional
differences can place their individual roughness coordinates on the same
reference scale, but their non-zero cross-order correlations remain and
produce the additional variance term
$2\sum_{r\neq s}\alpha_r\alpha_sG_{rs}^2$. This excess cannot be removed by
single-order normalization; it vanishes precisely under mutual
decorrelation. The operators of Section~\ref{sec:2} therefore combine three
properties: they retain their finite-difference order, have a common
unit-variance scale, and eliminate cross-order variance inflation under the
reference law.

\subsection{Sequential and simultaneous multi-order smoothing}

These properties motivate penalization across several difference orders. We
consider two distinct implementations: sequential smoothing, in which
order-specific smoothing operators are applied successively, and simultaneous
smoothing, in which all orders enter a single aggregate criterion. Both may
be defined using conventional, uncorrelated, or interpolated penalties.

For $r=1,\ldots,R$, let
\[
\mathfrak P_\eta^{(r)}
=
(1-\eta)D^{(r)\top}D^{(r)}
+
\eta\bar D^{(r)\top}\bar D^{(r)},
\qquad 0\leq\eta\leq1.
\]

\paragraph{Sequential smoothing.}
Starting from $f_{R+1}=X_i$, define successively, from order $R$ to order
$1$,
\[
f_r
=
\arg\min_{f\in\mathbb R^d}
\frac12
\left\{
\|f-f_{r+1}\|^2
+
\alpha_r f^\top\mathfrak P_\eta^{(r)}f
\right\}.
\]
The solution is
$f_r=(I_d+\alpha_r\mathfrak P_\eta^{(r)})^{-1}f_{r+1}$, and hence
\[
\hat X_i^{\mathrm{seq}}
=
\left[
\prod_{r=1}^{R}
\left(I_d+\alpha_r\mathfrak P_\eta^{(r)}\right)^{-1}
\right]X_i,
\]
where the product follows the order of the successive updates. The procedure
therefore applies an order-specific contraction at each stage.

\paragraph{Simultaneous smoothing.}
Alternatively, define the aggregate penalty
\[
\mathfrak P_{\boldsymbol\alpha,\eta}^{(1:R)}
=
\sum_{r=1}^R\alpha_r\mathfrak P_\eta^{(r)}.
\]
The simultaneous estimator is
\[
\hat X_i^{\mathrm{sim}}
=
\arg\min_{f\in\mathbb R^d}
\frac12
\left\{
\|X_i-f\|^2
+
f^\top\mathfrak P_{\boldsymbol\alpha,\eta}^{(1:R)}f
\right\}
=
\left(I_d+\mathfrak P_{\boldsymbol\alpha,\eta}^{(1:R)}\right)^{-1}X_i.
\]

The two constructions generally define different linear smoothers.
Sequential smoothing composes the order-specific contractions, whereas
simultaneous smoothing adds the corresponding quadratic penalties before
inversion. Even when the order-specific penalty matrices commute, the two
procedures need not coincide. If they share an eigenvector with eigenvalues
$\lambda_r$, their respective shrinkage factors along that direction are
\[
\prod_{r=1}^R(1+\alpha_r\lambda_r)^{-1}
\quad\text{and}\quad
\left(1+\sum_{r=1}^R\alpha_r\lambda_r\right)^{-1},
\]
respectively, which are generally different. Thus the two procedures should
be regarded as genuinely distinct regularization strategies rather than
alternative implementations of the same estimator.

When $\eta=0$, both constructions use conventional finite-difference
penalties, whereas $\eta=1$ gives the normalized and mutually decorrelated
operators studied above. The intermediate values $0<\eta<1$ provide a
computational interpolation between these endpoints; the exact reference
distribution in Proposition~\ref{prop:multiorder} and the variance-optimality
statement in Proposition~\ref{prop:variance-optimality} apply specifically
to the fully decorrelated endpoint $\eta=1$.

At this fully decorrelated endpoint, we write the simultaneous penalty and
smoother compactly as
\begin{equation}
\label{eq:simultaneous-Q0}
\mathfrak P_{\boldsymbol\alpha,Q_0}
:=
\mathfrak P_{\boldsymbol\alpha,1}^{(1:R)}
=
\sum_{r=1}^R\alpha_r\mathfrak P_{Q_0}^{(r)},
\qquad
S_{\boldsymbol\alpha,Q_0}
:=
\left(I_d+\mathfrak P_{\boldsymbol\alpha,Q_0}\right)^{-1},
\end{equation}
so that $\hat X_i^{\mathrm{sim}}=S_{\boldsymbol\alpha,Q_0}X_i$ at $\eta=1$.
This is the specific construction generalized to the covariance-adapted
reference geometry $Q_\vartheta$ in Section~\ref{sec:4}, by the substitution
$\bar D^{(r)}\to\bar D_\vartheta^{(r)}$.

For fixed $R$ and fixed stencil widths, all matrices
$\mathfrak P_\eta^{(r)}$ are banded, and their aggregate has bandwidth
bounded independently of $d$. The corresponding linear systems can therefore
be solved in $O(d)$ operations using standard banded-matrix methods. The
construction requires no basis expansion, knot placement, or
basis-dimension selection.

\section{Covariance-adapted statistical differentiation}
\label{sec:4}

The preceding construction is calibrated against the white-noise reference
law $Q_0=\mathcal N(0,I_d)$. When the reference noise exhibits dependence,
re-deriving a decorrelated stencil family from scratch for each covariance
structure by the recursive scheme of Section~\ref{sec:2} would be neither
systematic nor practical. We therefore allow the reference covariance
geometry to vary within a single parametric family
\[
Q_\vartheta=\mathcal N(0,\Sigma_\vartheta),
\qquad
\vartheta\in\Theta\subset\mathbb R^k,
\]
of which $Q_0$ is the special case $\Sigma_\vartheta=I_d$. We assume that
$\Theta$ is open, $\vartheta\mapsto\Sigma_\vartheta$ is continuously
differentiable, $\Sigma_\vartheta\succ0$ for every $\vartheta\in\Theta$,
$\operatorname{diag}(\Sigma_\vartheta)=\mathbf 1$ for every
$\vartheta\in\Theta$, and $\Sigma_\vartheta$ is translation-invariant on the
interior of the grid,
\[
(\Sigma_\vartheta)_{ij}=\gamma_\vartheta(i-j),
\]
for some covariance function $\gamma_\vartheta$ and every $i,j$ in the range
at which the stencils below are fully supported; thus $Q_\vartheta$ is a
standardized, translation-invariant covariance family. For a local stencil
$u$ (a finite sequence of coefficients indexed by offset $\ell$), let
$P_tu\in\mathbb R^d$ denote its zero-padded embedding at interior grid
location $t$, $(P_tu)_{t+\ell}=u_\ell$. Then, for local stencils $u,v$ and
every interior $t$,
\[
(P_tu)^\top\Sigma_\vartheta(P_tv)
=
\sum_{\ell,m}u_\ell v_m\,\gamma_\vartheta\{(t+\ell)-(t+m)\}
=
\sum_{\ell,m}u_\ell v_m\,\gamma_\vartheta(\ell-m),
\]
independent of $t$: this is exactly the property that allows a single
stencil, once constructed, to be validly slid across every interior grid
location. We therefore write $u^\top\Sigma_\vartheta v$ for this common
value, understanding $u,v$ to be embedded at any fixed interior location.
This induces the reference inner product
\[
\langle u,v\rangle_\vartheta
=
u^\top\Sigma_\vartheta v,
\]
which replaces the Euclidean geometry associated with $Q_0$.
An analogous replacement of Euclidean orthogonality by orthogonality with
respect to a covariance-induced geometry arises in robust sphering for
independent component analysis \citep{Vidal:2026}.

\subsection{Admissible difference spaces}

This section does not introduce a new differentiation scheme: it
keeps the admissible translated-difference spaces of Section~\ref{sec:2}
fixed and replaces Euclidean orthogonality by orthogonality under
$\langle u,v\rangle_\vartheta=u^\top\Sigma_\vartheta v$.

For a finite translation window $J_r\subset\mathbb Z$, let
\[
\mathcal V_r(J_r)
=
\operatorname{span}
\left\{
T^j\Delta^{(r)}:j\in J_r
\right\}
\]
denote the admissible space of translated $r$th-order differences introduced
in Proposition~\ref{prop:1}. We further define its parity subspace by
\[
\mathcal V_r^\pm(J_r)
=
\left\{
w\in\mathcal V_r(J_r):
w_{-\ell}=(-1)^r w_\ell
\text{ for all }\ell
\right\}.
\]
Write
$m_r=\dim\mathcal V_r^\pm(J_r)$ and let
$B_r\in\mathbb R^{d\times m_r}$ have columns forming a basis of this space.
Both $\mathcal V_r^\pm(J_r)$ and $B_r$ are independent of $\vartheta$.

With $w_\vartheta^{(0)}=e_0$ --- which, since
$\operatorname{diag}(\Sigma_\vartheta)=\mathbf 1$, already satisfies
$w_\vartheta^{(0)\top}\Sigma_\vartheta w_\vartheta^{(0)}=1$, consistent with
the normalization imposed at every subsequent order --- suppose that
$w_\vartheta^{(0)},\ldots,w_\vartheta^{(r-1)}$ have already been constructed.
We seek
$w_\vartheta^{(r)}=B_ra_r(\vartheta)$, where $a_r(\vartheta)$ is the
coordinate vector of the candidate stencil in the fixed basis $B_r$.
Its dependence on $\vartheta$ arises through the covariance-dependent
orthogonality conditions below. These coordinates should not be confused
with the penalty weights $\alpha_r$ of Section~\ref{sec:3}, nor with the
translated-difference coefficients $a_j^{(r)}(\vartheta)$ appearing in
Proposition~\ref{prop:cov-consistency}. The vector $a_r(\vartheta)$ must satisfy
\[
C_r(\vartheta)a_r(\vartheta)=0,
\qquad
C_r(\vartheta)
=
\begin{bmatrix}
w_\vartheta^{(0)\top}\Sigma_\vartheta B_r\\
\vdots\\
w_\vartheta^{(r-1)\top}\Sigma_\vartheta B_r
\end{bmatrix}
\in\mathbb R^{r\times m_r},
\]
together with the normalization
$w_\vartheta^{(r)\top}\Sigma_\vartheta w_\vartheta^{(r)}=1$ and a fixed sign
convention. Row $q+1$ of $C_r(\vartheta)a_r(\vartheta)=0$ reads
$w_\vartheta^{(q)\top}\Sigma_\vartheta w_\vartheta^{(r)}=0$, so the system is
exactly the $r$ orthogonality conditions against the previously constructed
orders $0,\ldots,r-1$. Thus the construction preserves the translated-difference class
while replacing Euclidean orthogonality by orthogonality with respect to the
reference covariance geometry.

\begin{proposition}
\label{prop:cov-consistency}
Whenever $w_\vartheta^{(r)}$ is defined by the construction,
$w_\vartheta^{(r)}\in\mathcal V_r(J_r)$, for every choice of $J_r$.
Consequently, wherever the construction is defined, the moment and
Taylor-expansion conclusions of Proposition~\ref{prop:1} hold independently
of the reference covariance used to determine the coefficients.
In particular,
\[
M_q(w_\vartheta^{(r)})=0,
\qquad q<r,
\]
and
\[
M_r(w_\vartheta^{(r)})
=
r!\sum_j a_j^{(r)}(\vartheta),
\]
whenever the latter quantity is nonzero. The parity condition also yields the
same order improvement described in Proposition~\ref{prop:1}.
\end{proposition}

The preceding result separates two aspects of the construction. The
finite-difference order is determined entirely by the admissible algebraic
space $\mathcal V_r(J_r)$, whereas the particular stencil selected within
that space is determined by the reference covariance geometry.

$\ker C_r(\vartheta)$ contains precisely the admissible coefficient vectors
satisfying all previous-order orthogonality constraints, so
$\dim\ker C_r(\vartheta)=1$ says exactly that there is a unique admissible
direction; normalization then fixes its magnitude and the sign convention
removes the remaining $\pm$ ambiguity. Since $\vartheta$ is later replaced
by an estimator $\widehat\vartheta_n$ (Section~\ref{sec:5}), we also require
the selected stencil to vary smoothly with $\vartheta$. The
one-dimensional-kernel condition is imposed throughout a neighbourhood
$U$ of $\tau$, rather than only at $\tau$, so that this unique direction
persists locally.

\begin{proposition}
\label{prop:cov-smoothness}
Fix translation windows $J_1,\ldots,J_R$ and let
$\tau\in\Theta$. Suppose that, for every $r=1,\ldots,R$,
\[
\dim\ker C_r(\vartheta)=1
\]
throughout a neighbourhood $U$ of $\tau$. Then
$w_\vartheta^{(r)}$ is uniquely determined up to sign and the mapping
$\vartheta\mapsto w_\vartheta^{(r)}$ is continuously differentiable on $U$
wherever the coordinate used to fix its sign does not vanish.
\end{proposition}

Whether $\dim\ker C_r(\vartheta)=1$ actually holds depends on the number
of effective constraints, rather than on the row count of $C_r(\vartheta)$
alone: some constraints may vanish identically for structural reasons,
leaving a larger kernel than the row count would suggest. The following
elementary fact identifies this source of redundancy under
reflection-symmetric reference geometries.

\begin{lemma}
\label{lem:crossparity}
Let $J$ denote reflection, $(Jx)_\ell=x_{-\ell}$, and suppose
\[
J\Sigma_\vartheta J=\Sigma_\vartheta.
\]
If $Ju=u$ and $Jv=-v$, then
$u^\top\Sigma_\vartheta v=0$. In particular, this applies to symmetric
Toeplitz covariance matrices.
\end{lemma}

\begin{remark}
\label{rem:cov-support}
Since $C_r(\vartheta)$ has $m_r$ columns and $r$ rows,
\[
\dim\ker C_r(\vartheta)\ge m_r-r.
\]
Thus the uniqueness condition of
Proposition~\ref{prop:cov-smoothness} necessarily requires
$m_r\le r+1$. Under a reflection-symmetric reference covariance,
Lemma~\ref{lem:crossparity} further implies that constraints involving
previous operators of opposite parity vanish identically. The effective
number of orthogonality constraints can therefore be strictly smaller than
$r$.

Consequently, enlarging the admissible parity space need not preserve
uniqueness and may instead increase the dimension of the solution set.
In particular, covariance orthogonality and derivative consistency alone need
not characterize a unique higher-order stencil once the admissible space
contains additional degrees of freedom. Thus the particular higher-order
stencils reported by Mizuta may require an additional selection criterion
beyond the orthogonality conditions considered here. Characterizing such a
criterion for non-minimal admissible spaces is left open.
\end{remark}

For $\vartheta$ satisfying the rank condition of
Proposition~\ref{prop:cov-smoothness} at order $r$, let
$\bar D_\vartheta^{(r)}\in\mathbb R^{d_r\times d}$ denote the
covariance-adapted difference matrix obtained by sliding the stencil
$w_\vartheta^{(r)}$ across the interior grid locations, in the same manner
as the operator $D^{(r)}$ of Section~\ref{sec:3} and the white-reference
matrix $\bar D^{(r)}$ of Section~\ref{sec:2}:
\[
(\bar D_\vartheta^{(r)}x)_t
=
\sum_\ell (w_\vartheta^{(r)})_\ell\, x_{t+\ell},
\qquad
t\in\mathcal I^{(r)}.
\]

Despite the possible non-uniqueness at higher orders described above,
whenever the condition of Proposition~\ref{prop:cov-smoothness} holds, the
decorrelation geometry of Section~\ref{sec:3} extends directly to the
covariance-adapted construction.

\begin{corollary}
\label{cor:cov-geometry}
Suppose the conditions of Proposition~\ref{prop:cov-smoothness} hold and let
$\varepsilon\sim Q_\vartheta$ denote the corresponding reference vector, as in
Section~\ref{sec:2}. At every grid location $t$ at which all stencils are
fully supported,
\[
Z_\vartheta(t)
=
\begin{pmatrix}
(\bar D_\vartheta^{(1)}\varepsilon)_t\\
\vdots\\
(\bar D_\vartheta^{(R)}\varepsilon)_t
\end{pmatrix}
\sim\mathcal N(0,I_R).
\]
\end{corollary}

\subsection{Covariance-adapted penalization and smoothing}

We now connect this construction to the discrete penalization framework of
Section~\ref{sec:3} by the substitution $\bar D^{(r)}\to\bar D_\vartheta^{(r)}$.
Exactly as $\mathfrak P_{Q_0}^{(r)}=\bar D^{(r)\top}\bar D^{(r)}$ and, at the
fully decorrelated endpoint $\eta=1$,
$\mathfrak P_{\boldsymbol\alpha,Q_0}=\sum_{r=1}^R\alpha_r\mathfrak P_{Q_0}^{(r)}$
and $S_{\boldsymbol\alpha,Q_0}=(I_d+\mathfrak P_{\boldsymbol\alpha,Q_0})^{-1}$
define the white-reference simultaneous penalty and smoother
(Section~\ref{sec:3}), the covariance-adapted versions are, for $\vartheta$
satisfying the rank condition of Proposition~\ref{prop:cov-smoothness} at
every order up to $R$ and fixed non-negative weights
$\boldsymbol\alpha=(\alpha_1,\ldots,\alpha_R)^\top$,
\begin{equation}
\label{eq:cov-adapted-smoother}
\mathfrak P_{\boldsymbol\alpha,\vartheta}
=
\sum_{r=1}^R
\alpha_r
\bar D_\vartheta^{(r)\top}
\bar D_\vartheta^{(r)},
\qquad
S_{\boldsymbol\alpha,\vartheta}
=
\left(I_d+\mathfrak P_{\boldsymbol\alpha,\vartheta}\right)^{-1}.
\end{equation}
The weights $\boldsymbol\alpha$ are held fixed throughout the remainder of
the paper; from here on we suppress them from the notation and write
\[
\mathfrak P_\vartheta := \mathfrak P_{\boldsymbol\alpha,\vartheta},
\qquad
S_\vartheta := S_{\boldsymbol\alpha,\vartheta},
\qquad
\widehat f_\vartheta(x) := S_\vartheta x.
\]
In Section~\ref{sec:5}, $S_\vartheta$ will be evaluated at the oracle
parameter $\vartheta_0$ and at its estimator $\widehat\vartheta_n$;
relating the two through a first-order expansion requires the dependence
of the smoother on $\vartheta$ to be differentiable, not merely well
defined.

\begin{proposition}
\label{prop:S-smooth}
Suppose that the conditions of Proposition~\ref{prop:cov-smoothness} hold
throughout a neighbourhood $U$ of $\tau$ at every order $r=1,\ldots,R$.
Then $\vartheta\mapsto S_\vartheta$ is continuously differentiable on $U$,
with
\begin{equation}
\label{eq:Sderivative}
\dot S_{\vartheta,j}
=
-S_\vartheta\dot{\mathfrak P}_{\vartheta,j}S_\vartheta,
\qquad
\dot{\mathfrak P}_{\vartheta,j}
=
\sum_{r=1}^R\alpha_r
\left\{
\dot{\bar D}_{\vartheta,j}^{(r)\top}\bar D_\vartheta^{(r)}
+
\bar D_\vartheta^{(r)\top}\dot{\bar D}_{\vartheta,j}^{(r)}
\right\},
\qquad
\vartheta\in U.
\end{equation}
\end{proposition}

\section{Estimation of the reference geometry}
\label{sec:5}

The covariance-adapted operators of Section~\ref{sec:4} depend on the
reference parameter $\vartheta$. If $\vartheta_0$ were known, these
operators, and the corresponding smoother $S_{\vartheta_0}$, could be
constructed directly; in practice $\vartheta_0$ is typically unknown and
must instead be estimated before they can be implemented. Throughout this
section, the grid size $d$, the number of orders $R$, the stencil
translation windows $J_1,\ldots,J_R$, and the regularization weights
$\boldsymbol\alpha$ are all fixed, as in Sections~\ref{sec:2}--\ref{sec:4};
the only asymptotic index is the calibration sample size $n\to\infty$, and
$n$ refers exclusively to the number of independent draws
$\varepsilon_1,\ldots,\varepsilon_n\sim Q_{\vartheta_0}$ used to estimate
$\vartheta_0$, not to the observed curve $X$ (or $X_i$) being regularized,
which remains a single fixed object throughout.
We first study estimation of the reference geometry through the
multi-order roughness coordinates themselves and then quantify the
first-order effect of replacing the oracle parameter by its estimate in the
resulting curve regularization procedure.

\subsection{Estimating equations for the reference geometry}

Fix an interior grid location $t$ at which all stencils
$w_\vartheta^{(1)},\ldots,w_\vartheta^{(R)}$ are fully supported, and
identify each $w_\vartheta^{(r)}$ with its zero-padded vector in
$\mathbb R^d$ centred at $t$, as in Corollary~\ref{cor:cov-geometry}. For
$x\in\mathbb R^d$, define
\[
h_\vartheta(x)
=
\big(
(w_\vartheta^{(1)\top}x)^2-1,\ldots,
(w_\vartheta^{(R)\top}x)^2-1
\big)^\top
\in\mathbb R^R,
\qquad
\mathcal R_{\boldsymbol\alpha}(x;\vartheta,t)
=
\sum_{r=1}^R
\alpha_r
(w_\vartheta^{(r)\top}x)^2,
\]
so that $\boldsymbol\alpha^\top h_\vartheta(x)=\mathcal
R_{\boldsymbol\alpha}(x;\vartheta,t)-\sum_r\alpha_r$; this extends the
white-reference roughness statistic $\mathcal R_{\boldsymbol\alpha}(\cdot;t)$
of Section~\ref{sec:3} to the covariance-adapted geometry $Q_\vartheta$. If
$\varepsilon\sim Q_\vartheta$, then by construction
$E_\vartheta\{h_\vartheta(\varepsilon)\}=0$ for every
$\vartheta\in\Theta$. Thus the order-specific standardized roughness
coordinates, evaluated at the reference law, provide a family of estimating
equations for the reference covariance parameter. Their informativeness about
$\vartheta$ is determined by the sensitivity matrix introduced below.

\begin{proposition}
\label{prop:Dtau}
Suppose that $\vartheta\mapsto w_\vartheta^{(r)}$ is continuously
differentiable at $\tau$ for every $r=1,\ldots,R$. Then
\begin{equation}
\label{eq:Dtau}
\Gamma_{\tau,rj}
:=
-E_\tau
\left\{
h_{\tau,r}(\varepsilon)\kappa_{\tau,j}(\varepsilon)
\right\}
=
-w_\tau^{(r)\top}
\dot\Sigma_{\tau,j}
w_\tau^{(r)},
\end{equation}
where
$\dot\Sigma_{\tau,j}
=
\partial_{\vartheta_j}\Sigma_\vartheta|_{\vartheta=\tau}$
and $\kappa_\tau$ denotes the score of
$\mathcal N(0,\Sigma_\vartheta)$ at $\tau$.
\end{proposition}

Thus the sensitivity of the order-$r$ roughness condition to
$\vartheta_j$ is determined directly by the variation of the reference
covariance geometry along the direction $w_\tau^{(r)}$.

The covariance structure of these estimating equations has an equally
explicit form.

\begin{proposition}
\label{prop:Omega}
Let $G_\tau$ denote the $\Sigma_\tau$-correlation matrix of the normalized
stencils,
\[
(G_\tau)_{rs}
=
w_\tau^{(r)\top}\Sigma_\tau w_\tau^{(s)}.
\]
Then
\begin{equation}
\label{eq:Omega}
\Omega_\tau
:=
\operatorname{cov}_\tau\{h_\tau(\varepsilon)\}
=
2(G_\tau\circ G_\tau),
\end{equation}
where $\circ$ denotes the Hadamard product. In particular, mutual
decorrelation gives
$\Omega_\tau=2I_R$.
\end{proposition}

The expression \eqref{eq:Omega} uses Gaussian fourth-moment identities.
For non-Gaussian reference laws, decorrelation of the linear roughness
coordinates need not imply $\Omega_\tau=2I_R$. Since $G_\tau$ is the Gram
matrix of $w_\tau^{(1)},\ldots,w_\tau^{(R)}$ under the inner product
$\langle u,v\rangle_\tau=u^\top\Sigma_\tau v$, it is positive definite
whenever these stencils are linearly independent. In particular, the
moment conditions of Proposition~\ref{prop:cov-consistency} imply linear
independence whenever $M_r(w_\tau^{(r)})\ne0$ for every $r$: if $r_0$ is
the smallest index with nonzero coefficient in a linear dependence,
applying $M_{r_0}$ leaves only the $r_0$th term, a contradiction. Hence,
by the Schur product theorem, $\Omega_\tau=2(G_\tau\circ G_\tau)\succ0$,
so $\Omega_\tau^{-1}$ is well defined wherever it is used below.

\begin{remark}
\label{rem:two-consequences}
Corollary~\ref{cor:cov-geometry} and
Proposition~\ref{prop:Omega} describe two distinct consequences of
\[
w_\vartheta^{(r)\top}
\Sigma_\vartheta
w_\vartheta^{(s)}
=
\delta_{rs}.
\]
The former diagonalizes the local roughness geometry, giving
$Z_\vartheta(t)\sim\mathcal N(0,I_R)$, whereas the latter diagonalizes the
covariance of the estimating equations used to calibrate that geometry,
giving $\Omega_\vartheta=2I_R$. Both are exact finite-dimensional
properties.
\end{remark}

When the number of roughness conditions exceeds the dimension of
$\vartheta$, this covariance structure has a direct consequence for their
combination.

\begin{corollary}
\label{cor:efficiency}
Suppose $R>k$, $\Gamma_\tau$ has full column rank $k$, and the remaining
regularity conditions for GMM estimation hold. Among GMM estimators based on
$h_\vartheta$ with weight matrix $\mathcal W\in\mathbb R^{R\times R}$,
$\mathcal W\succ0$, the efficient choice is
$\mathcal W=\Omega_\tau^{-1}$, up to multiplication by a positive scalar
\citep[Thm.~3.2]{Hansen:1982}.
Under mutual decorrelation, $\Omega_\tau=2I_R$, and hence the unweighted
choice $\mathcal W=I_R$ already attains the efficient GMM covariance within
this moment family.
\end{corollary}

The efficiency statement in Corollary~\ref{cor:efficiency} is relative to
estimators constructed from these particular roughness conditions; no claim
of Fisher or semiparametric efficiency in a larger statistical model is
intended.

\subsection{Plug-in estimation of the reference geometry}

Let $\varepsilon_1,\ldots,\varepsilon_n$ be an independent calibration sample
from $Q_{\vartheta_0}$, distinct from the observed curves $X_i$ to be
regularized, and define
\[
H_n(\vartheta)
=
\frac{1}{\sqrt n}
\sum_{i=1}^n h_\vartheta(\varepsilon_i).
\]
Let $\widehat\vartheta_n$ be an estimator of $\vartheta_0$ constructed from
these estimating equations.

\begin{remark}
\label{rem:consistency}
Throughout the remainder of this section, we assume that the GMM
estimator $\widehat\vartheta_n$ is consistent and
satisfies
\[
\widehat\vartheta_n-\vartheta_0=O_p(n^{-1/2}).
\]
These properties hold under the usual identification, smoothness, and
uniform-convergence conditions for finite-dimensional GMM estimation,
applied to the sample moment $n^{-1}\sum_{i=1}^n h_\vartheta(\varepsilon_i)$;
we do not re-derive them here. The corresponding GMM theorem is invoked
below only for its stated purpose, the asymptotic normal distribution of
$\widehat\vartheta_n$ (Corollary~\ref{cor:vartheta-normal}).
\end{remark}

\begin{remark}
\label{rem:qmd-automatic}
Because $\Theta$ is open, $\Sigma_\vartheta\succ0$, and
$\vartheta\mapsto\Sigma_\vartheta$ is continuously differentiable
(Section~\ref{sec:4}), the Gaussian family
$\{\mathcal N(0,\Sigma_\vartheta):\vartheta\in\Theta\}$ is differentiable
in quadratic mean at every $\vartheta_0\in\Theta$. Indeed, the density is
everywhere differentiable in $\vartheta$ with score $\kappa_{\vartheta,j}$
as in the proof of Proposition~\ref{prop:Dtau}, and the resulting Fisher
information, whose $(j,k)$ entry is
\[
\tfrac12\operatorname{tr}\!\left(
\Sigma_\vartheta^{-1}\dot\Sigma_{\vartheta,j}
\Sigma_\vartheta^{-1}\dot\Sigma_{\vartheta,k}
\right),
\]
is continuous on $\Theta$ under the stated $C^1$ and positive-definiteness
assumptions; this is precisely the
standard sufficient condition for differentiability in quadratic mean
\citep[Lemma~7.6]{vanderVaart:1998}. In the terminology of
\citet{Schick:2001}, this is precisely Hellinger differentiability, which
is the property required below. Hellinger differentiability is therefore
not an independent hypothesis in what follows, but a standing consequence
of the assumptions already imposed on $\Sigma_\vartheta$.
\end{remark}

\begin{theorem}
\label{thm:plugin}
Suppose that the rank condition of
Proposition~\ref{prop:cov-smoothness} holds throughout a neighbourhood of
$\vartheta_0$ at every order. Then, since
$\{\mathcal N(0,\Sigma_\vartheta):\vartheta\in\Theta\}$ is Hellinger
differentiable at $\vartheta_0$ (Remark~\ref{rem:qmd-automatic}),
$\vartheta\mapsto h_\vartheta$ is Hellinger continuous at
$\vartheta_0$, and
\begin{equation}
\label{eq:plugin-H}
H_n(\widehat\vartheta_n)
-
H_n(\vartheta_0)
-
\Gamma_{\vartheta_0}
\sqrt n(\widehat\vartheta_n-\vartheta_0)
\xrightarrow{P_{\vartheta_0}}0.
\end{equation}
\end{theorem}

For
\[
\bar{\mathcal R}_{\boldsymbol\alpha}(\vartheta;t)
=
\frac1n\sum_{i=1}^n
\mathcal R_{\boldsymbol\alpha}(\varepsilon_i;\vartheta,t),
\]
we have
$\boldsymbol\alpha^\top H_n(\vartheta)
=
\sqrt n\{\bar{\mathcal R}_{\boldsymbol\alpha}(\vartheta;t)
-\sum_r\alpha_r\}$.
Theorem~\ref{thm:plugin} therefore immediately yields the following
aggregate version.

\begin{corollary}
\label{cor:plugin-aggregate}
Under the conditions of Theorem~\ref{thm:plugin},
\begin{equation}
\label{eq:plugin-aggregate}
\sqrt n
\left\{
\bar{\mathcal R}_{\boldsymbol\alpha}
(\widehat\vartheta_n;t)
-
\bar{\mathcal R}_{\boldsymbol\alpha}
(\vartheta_0;t)
\right\}
-
\boldsymbol\alpha^\top \Gamma_{\vartheta_0}
\sqrt n(\widehat\vartheta_n-\vartheta_0)
\xrightarrow{P_{\vartheta_0}}0.
\end{equation}
\end{corollary}

Thus estimation of the reference geometry contributes a generally
non-negligible first-order term to the limiting distribution of plug-in
roughness statistics. In particular, this contribution must be incorporated
when such statistics are used for calibration or specification assessment.

\begin{corollary}
\label{cor:vartheta-normal}
Let $\widehat\vartheta_n$ denote the root-$n$-consistent estimator of
$\vartheta_0$ introduced above (Remark~\ref{rem:consistency}), obtained as a GMM
estimator based on the moment conditions $h_\vartheta$ with a fixed weight
matrix $\mathcal W\in\mathbb R^{R\times R}$, $\mathcal W\succ0$, and suppose
that $\Gamma_{\vartheta_0}$ has full column rank $k$ (cf.\
Corollary~\ref{cor:efficiency}). Then, under the standard identification and
regularity conditions for asymptotically normal GMM estimation
\citep[Thm.~3.1]{Hansen:1982}, specialized here to an independent and
identically distributed calibration sample and to the quadratic weighting
$\mathcal W$ used throughout this paper,
\begin{equation}
\label{eq:vartheta-normal}
\sqrt n(\widehat\vartheta_n-\vartheta_0)
\Rightarrow
\mathcal N\{0,V_{\vartheta_0}(\mathcal W)\},
\qquad
V_{\vartheta_0}(\mathcal W)
=
(\Gamma_{\vartheta_0}^\top \mathcal W\Gamma_{\vartheta_0})^{-1}
\Gamma_{\vartheta_0}^\top \mathcal W\Omega_{\vartheta_0}\mathcal W\Gamma_{\vartheta_0}
(\Gamma_{\vartheta_0}^\top \mathcal W\Gamma_{\vartheta_0})^{-1}.
\end{equation}
\end{corollary}

The expansion of Theorem~\ref{thm:plugin} and the asymptotic distribution of
Corollary~\ref{cor:vartheta-normal} are parallel consequences of the same
root-$n$-consistent estimator $\widehat\vartheta_n$, obtained by different
arguments for different purposes. Theorem~\ref{thm:plugin}, via the
framework of \citet{Schick:2001}, controls the empirical estimating
equations $H_n$ evaluated \emph{at} $\widehat\vartheta_n$ rather than at
$\vartheta_0$, as required for the roughness aggregate of
Corollary~\ref{cor:plugin-aggregate}; Corollary~\ref{cor:vartheta-normal},
by contrast, supplies the limiting distribution of $\widehat\vartheta_n$
itself, from standard GMM theory.

\subsection{Covariance-adapted curve regularization}
\label{sec:5.3}

We now propagate estimation of the reference geometry to the curve
regularization procedure itself, using the covariance-adapted smoother
$S_\vartheta$ already defined in Section~\ref{sec:4}. Replacing the oracle
reference geometry $\vartheta_0$ by its estimate $\widehat\vartheta_n$
yields the plug-in regularized curve $\widehat X_i=S_{\widehat\vartheta_n}X_i$,
applied to a curve $X_i$ independent of the calibration sample, in place of
its oracle counterpart $S_{\vartheta_0}X_i$.

Recall that Proposition~\ref{prop:S-smooth} is an ordinary
finite-dimensional differentiability result, requiring only the $C^1$
dependence of $\Sigma_\vartheta$ established in Section~\ref{sec:4}.
Hellinger differentiability instead enters through Theorem~\ref{thm:plugin},
where it underlies the Schick-type linearization of the parameter-dependent
empirical estimating equations $H_n(\vartheta)$; propagating the resulting
behaviour of $\widehat\vartheta_n$ into $S_{\widehat\vartheta_n}$, by
contrast, uses only the ordinary matrix calculus of
Proposition~\ref{prop:S-smooth}.

\begin{proposition}
\label{prop:fhat-expansion}
Let $\widehat\vartheta_n$ satisfy
$\widehat\vartheta_n-\vartheta_0=O_p(n^{-1/2})$, and let $X$ be a new curve,
independent of the calibration sample $\varepsilon_1,\ldots,\varepsilon_n$
used to construct $\widehat\vartheta_n$, such that $\|X\|=O_p(1)$. Then
\begin{equation}
\label{eq:fhat-expansion}
\widehat f_{\widehat\vartheta_n}(X)
-
\widehat f_{\vartheta_0}(X)
=
-S_{\vartheta_0}
\left[
\sum_{j=1}^k
\dot{\mathfrak P}_{\vartheta_0,j}
(\widehat\vartheta_{n,j}-\vartheta_{0,j})
\right]
S_{\vartheta_0}X
+
o_p(n^{-1/2}).
\end{equation}
\end{proposition}

The correction in \eqref{eq:fhat-expansion} is generically of order
$n^{-1/2}$; in particular, no general property of the construction forces
$\dot{\mathfrak P}_{\vartheta_0,j}$ to vanish. Estimating the reference
geometry is therefore not first-order negligible. The relevant question is
instead how the leading-order cost of this estimation can be minimized.

For a fixed $x\in\mathbb R^d$, define
$L(x)\in\mathbb R^{d\times k}$ by
\begin{equation}
\label{eq:Lx}
L_{\cdot j}(x)
=
-S_{\vartheta_0}
\dot{\mathfrak P}_{\vartheta_0,j}
S_{\vartheta_0}x,
\qquad j=1,\ldots,k.
\end{equation}
By \eqref{eq:Sderivative},
$L_{\cdot j}(x)=\dot S_{\vartheta_0,j}x$, so that $L(x)$ is the
Jacobian of the map $\vartheta\mapsto S_\vartheta x$ at
$\vartheta_0$. Thus, $L(x)$ maps a local perturbation of the reference
covariance parameter to the corresponding first-order perturbation of the
regularized curve.

\begin{theorem}
\label{thm:efficient-plugin-risk}
Suppose $\widehat\vartheta_n$ and $\mathcal W$ satisfy the conditions of
Corollary~\ref{cor:vartheta-normal}, so that
$\sqrt n(\widehat\vartheta_n-\vartheta_0)\Rightarrow\mathcal
N\{0,V_{\vartheta_0}(\mathcal W)\}$, and suppose in addition that the sequence
$n\|\widehat f_{\widehat\vartheta_n}(x)-\widehat f_{\vartheta_0}(x)\|^2$,
$n\ge1$, is uniformly integrable. Then, for every fixed new curve
$x$,
\begin{equation}
\label{eq:conditional-plugin-risk}
nE
\left\|
\widehat f_{\widehat\vartheta_n}(x)
-
\widehat f_{\vartheta_0}(x)
\right\|^2
\longrightarrow
\operatorname{tr}
\left\{
L(x)^\top L(x)V_{\vartheta_0}(\mathcal W)
\right\}.
\end{equation}
Among GMM estimators based on the moment conditions $h_\vartheta$, the
right-hand side of \eqref{eq:conditional-plugin-risk} is minimized by the
efficient GMM weight
$\mathcal W=\Omega_{\vartheta_0}^{-1}$, up to multiplication by a positive scalar,
simultaneously for every $x$ and every fixed
$\boldsymbol\alpha$.

Under mutual decorrelation,
$\Omega_{\vartheta_0}=2I_R$, so the ordinary unweighted choice $\mathcal W=I_R$
already attains this minimum.
\end{theorem}

The preceding result is conditional on the curve being regularized. It also
immediately yields an integrated version when the new curve is random.

\begin{corollary}
\label{cor:integrated-plugin-risk}
Let $X$ be independent of the calibration sample
$\varepsilon_1,\ldots,\varepsilon_n$, with $\|X\|=O_p(1)$, suppose
\[
M_{\vartheta_0}
:=
E\{L(X)^\top L(X)\}
\]
is finite, and suppose that the sequence
$n\|\widehat f_{\widehat\vartheta_n}(X)-\widehat f_{\vartheta_0}(X)\|^2$,
$n\ge1$, is uniformly integrable (jointly over the randomness in $X$ and in
the calibration sample). Then
\begin{equation}
\label{eq:integrated-plugin-risk}
nE
\left\|
\widehat f_{\widehat\vartheta_n}(X)
-
\widehat f_{\vartheta_0}(X)
\right\|^2
\longrightarrow
\operatorname{tr}
\left\{
M_{\vartheta_0}V_{\vartheta_0}(\mathcal W)
\right\}.
\end{equation}
Consequently, the efficient GMM weight also minimizes the integrated
leading-order plug-in risk. Under mutual decorrelation this minimum is again
attained by the unweighted estimator.
\end{corollary}

\begin{remark}
\label{rem:plugin-risk-scope}
Theorem~\ref{thm:efficient-plugin-risk} and
Corollary~\ref{cor:integrated-plugin-risk} concern the additional cost of
estimating the reference geometry relative to the oracle covariance-adapted
smoother,
\[
\widehat f_{\widehat\vartheta_n}
-
\widehat f_{\vartheta_0}.
\]
They do not establish optimality of the decorrelated stencil family over
alternative families of difference operators, since changing the stencil
family changes the sensitivity matrix, the covariance of the estimating
equations, and the downstream smoother simultaneously. Nor do these results
characterize the total recovery risk
$E\|\widehat f_{\widehat\vartheta_n}-f\|^2$.
Rather, conditional on the chosen moment family, efficient estimation of the
reference geometry minimizes its leading-order propagation into the
regularized curve. For the decorrelated construction, the additional
property $\Omega_{\vartheta_0}=2I_R$ makes this efficient calibration
available through the ordinary unweighted estimating equations.
\end{remark}

\section{Selection of $\alpha$ by GCV}
\label{sec:6}

For each observation $X_i\in\mathbb R^d$ and a fixed penalty operator
$\mathfrak P$, the regularization parameter $\alpha$ can be selected by
minimizing a generalized cross-validation (GCV) criterion associated with the
discrete linear smoother
\[
\hat X_i(\alpha)=S_\alpha X_i,
\qquad
S_\alpha=(I_d+\alpha\mathfrak P)^{-1}.
\]
Here $\mathfrak P$ may represent either a single-order penalty or a
prespecified multi-order aggregate penalty. In the latter case, the relative
contributions of the different derivative orders are regarded as part of the
fixed penalty structure, while $\alpha$ controls its overall regularization
strength.

The GCV score is defined as
\[
\mathrm{GCV}_i(\alpha)
=
\frac{\|(I_d-S_\alpha)X_i\|^2}
     {(d-\operatorname{tr}(S_\alpha))^2},
\qquad
\hat\alpha_i\in
\arg\min_{\alpha\in\mathcal A}\mathrm{GCV}_i(\alpha),
\]
where $\mathcal A$ is a finite evaluation grid. The quantity
$\operatorname{tr}(S_\alpha)$ serves as an effective degrees-of-freedom
measure, reflecting the contraction induced by the smoother.

More generally, when the penalty belongs to a prescribed family
$\{\mathfrak P_\eta:\eta\in[0,1]\}$, the same criterion may be applied
conditionally on a fixed value of $\eta$, with
\[
S_{\alpha,\eta}
=
(I_d+\alpha\mathfrak P_\eta)^{-1}.
\]
Thus, the structural form of the penalty is specified separately, whereas GCV
is used to tune its overall regularization strength $\alpha$.

In the sequential scheme of Section~\ref{sec:3}, GCV may be applied
separately at each order, with $\mathfrak P_\eta^{(r)}$ fixed and the current
input $f_{r+1}$ replacing $X_i$ in the criterion, thereby yielding an
order-specific value $\hat\alpha_{i,r}$. In the simultaneous scheme, one may
analogously optimize over the full vector $\boldsymbol\alpha$. For
computational convenience, however, a simpler alternative is to prescribe
the relative multi-order penalty structure and use GCV only to select a
single overall regularization strength. For example, restricting the
order-specific weights to $\alpha_r=\alpha c_r$, with fixed non-negative
weights $c_1,\ldots,c_R$, gives
\[
\mathfrak P_\eta
=
\sum_{r=1}^R c_r\mathfrak P_\eta^{(r)},
\qquad
S_{\alpha,\eta}
=
\left(
I_d+\alpha\sum_{r=1}^R c_r\mathfrak P_\eta^{(r)}
\right)^{-1},
\]
so that the $c_r$ determine the relative contribution of the different
orders and $\alpha$ determines their common overall strength.

The numerical experiments below use oracle tuning when the target signal is
available by design, so that comparisons isolate reconstruction performance
from tuning performance rather than the performance of the tuning
criterion.

Here GCV is used operationally as a data-adaptive criterion for balancing
residual energy against smoothing strength within the finite-dimensional
observation space; no asymptotic optimality claim for the GCV-selected
regularization parameter is made.

\section{Numerical experiments}
\label{sec:7}

We use three numerical experiments to examine complementary aspects of the
theory. The first isolates the effect of the reference covariance geometry,
the second evaluates the cost of estimating that geometry, and the third
assesses practical curve reconstruction. An additional reconstruction
experiment for a smooth sinusoidal target is
reported in Appendix~\ref{appB}.

\paragraph{\textbf{Covariance geometry.}}
We first isolate the effect of covariance adaptation from that of tuning.
We consider an AR(1) covariance model on a grid of size $d=100$,
\[
\Sigma_\rho=(\rho^{|j-k|})_{j,k=1}^d,
\qquad
\rho\in\{0,0.3,0.6,0.9\},
\]
and construct the covariance-adapted difference operators
$\bar D_\rho^{(r)}$, $r=1,2,3$, as described in Section~\ref{sec:4}.
The same regularization parameters are used throughout the experiment:
they are calibrated once at $\rho=0.6$ by oracle minimization of the
reconstruction MSE against the known target and subsequently held fixed for
all values of $\rho$. Hence, changes in the resulting smoothing operators are
due solely to changes in covariance geometry, rather than to retuning of
the regularization parameters.

For each value of $\rho$, we compare the oracle smoother
$S_\rho$, constructed using the true covariance parameter, with two
alternatives. The first is the independence-based smoother $S_0$, which
ignores the correlation structure. The second is the plug-in smoother
$S_{\widehat\rho}$, where $\widehat\rho$ is estimated from an independent
calibration sample $\varepsilon_1,\ldots,\varepsilon_{200}\sim Q_\rho$
using the moment criterion of Section~\ref{sec:5}. Results are averaged
over 100 Monte Carlo replications.

To measure the effect at the operator level, we report the normalized
squared Frobenius discrepancies
\[
d^{-1}\|S_0-S_\rho\|_{\mathrm F}^2
\qquad\text{and}\qquad
d^{-1}\|S_{\widehat\rho}-S_\rho\|_{\mathrm F}^2.
\]
The former measures the distortion induced by ignoring the covariance
geometry, whereas the latter measures the additional error introduced by
estimating it. At $\rho=0$, the oracle and independence-based constructions
coincide exactly. As the dependence increases, the discrepancy between
$S_0$ and $S_\rho$ increases sharply, while the plug-in smoother remains
close to its oracle counterpart; see \autoref{tab:geometry}.

\begin{table}[t]
\caption{Covariance-geometry experiment. RMSE of the estimated AR(1)
parameter and base-10 logarithms of the normalized squared Frobenius
discrepancies of the independence-based and plug-in operators from the
oracle smoothing operator. Results are averaged over 100 replications.}
\label{tab:geometry}
\centering
\begin{tabular}{@{}cccc@{}}
\hline
$\rho$
& $\operatorname{RMSE}(\widehat\rho)$
& $\log_{10}\!\left\{d^{-1}\|S_0-S_\rho\|_{\mathrm F}^2\right\}$
& $\log_{10}\!\left\{d^{-1}\|S_{\widehat\rho}-S_\rho\|_{\mathrm F}^2\right\}$ \\
\hline
0.0 & 0.0411 & ---     & $-7.27$ \\
0.3 & 0.0392 & $-5.30$ & $-6.73$ \\
0.6 & 0.0241 & $-4.18$ & $-6.17$ \\
0.9 & 0.0065 & $-2.65$ & $-5.01$ \\
\hline
\end{tabular}

\smallskip
\footnotesize
At $\rho=0$, $S_0=S_\rho$ exactly, so the corresponding discrepancy is
zero and its logarithm is not reported. Since $\rho=0$ is a boundary point
of the nonnegative AR(1) parameter space, its RMSE is reported only as a
finite-sample diagnostic and is not used to assess the interior-point
asymptotic theory.
\end{table}

\paragraph{\textbf{Estimated versus oracle geometry.}}
We next examine the effect of estimating the covariance geometry, as opposed
to treating it as known. We fix an AR(1) covariance model with parameter
$\rho_0=0.6$ on a grid of size $d=100$ and construct the oracle smoother
$S_{\rho_0}$ using the true covariance parameter. For each calibration
sample size $n$, an independent calibration sample
$\varepsilon_1,\ldots,\varepsilon_n\sim Q_{\rho_0}$ is used to estimate
the AR(1) parameter, yielding $\widehat\rho_n$ and the plug-in smoother
$S_{\widehat\rho_n}$.

To isolate the error due exclusively to covariance estimation, the
regularization parameters are held fixed and both smoothers are applied to
the same independently generated curve
\[
Y=f+0.2\,\varepsilon_{\mathrm{new}},
\qquad
f=(\sin(6\pi x_t))_{t\in\mathcal I_d},
\]
where $(x_t)_{t\in\mathcal I_d}$ is an equally spaced grid on $[0,1]$ and
$\varepsilon_{\mathrm{new}}$ has AR(1) covariance with parameter $\rho_0$
and is independent of the calibration sample. We then compute
\[
\Delta_n
=
\frac{1}{d}
\left\|
S_{\widehat\rho_n}Y-S_{\rho_0}Y
\right\|^2.
\]
Thus, unlike reconstruction error, $\Delta_n$ vanishes when the estimated
and oracle covariance geometries coincide and does not confound
geometry-estimation error with the smoothing error already present in the
oracle reconstruction $S_{\rho_0}Y$.

\begin{figure}[t]
\centering
\includegraphics[width=\textwidth]{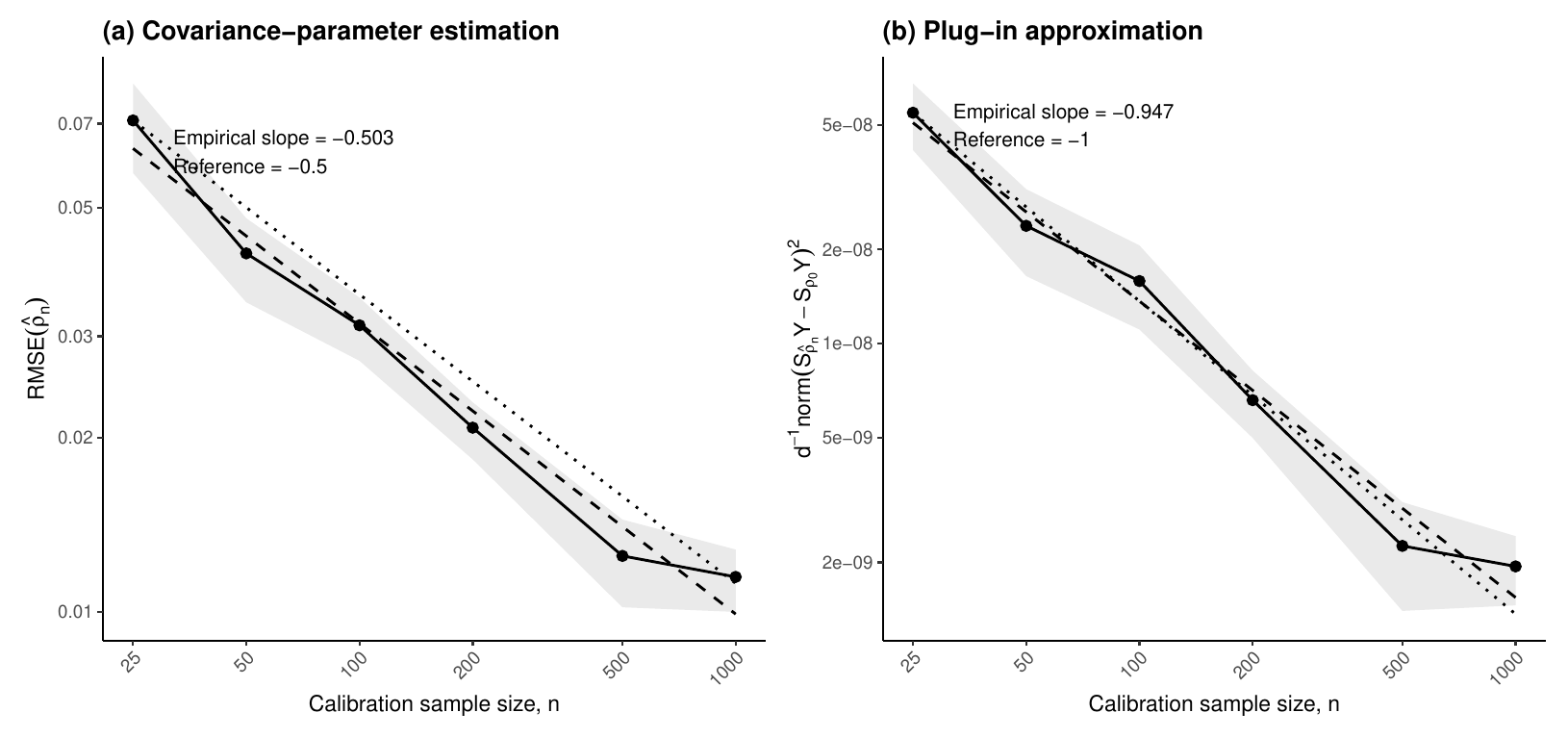}
\caption{Estimated versus oracle covariance geometry under an AR(1) model
with $\rho_0=0.6$.
Panel (a) shows the Monte Carlo RMSE of the covariance-parameter estimator
$\widehat\rho_n$ as a function of the calibration sample size $n$.
(b) Mean squared discrepancy
$d^{-1}\|S_{\widehat\rho_n}Y-S_{\rho_0}Y\|^2$
between the plug-in and oracle smoothers.
Both axes are logarithmic. Solid lines and points denote the
estimates, with 95\% Monte Carlo intervals shown as shaded bands; dashed
lines are fitted log--log regressions, and dotted lines show the
theoretical $n^{-1/2}$ and $n^{-1}$ reference rates, respectively.
The fitted slopes are $-0.503$ in (a) and $-0.947$ in (b).}
\label{fig:plugin_geometry}
\end{figure}

According to the plug-in expansion of Section~\ref{sec:5},
$\widehat\rho_n-\rho_0=O_p(n^{-1/2})$ implies
$\Delta_n=O_p(n^{-1})$. For each calibration sample size, we use 100 Monte
Carlo replications. We assess the predicted rates by fitting log--log
regressions to the RMSE of $\widehat\rho_n$ and to $\Delta_n$ across $n$.
The resulting empirical slopes are $-0.503$ and $-0.947$, respectively,
closely matching the theoretical reference slopes $-1/2$ and $-1$; see
\autoref{fig:plugin_geometry}. This provides a direct numerical check of
both the root-$n$ estimation of the covariance parameter and the resulting
first-order plug-in approximation.

\paragraph{\textbf{Reconstruction performance.}}
Finally, we examine reconstruction performance under correlated noise and
compare the covariance-adapted construction with conventional discrete,
basis, and kernel smoothers. The target curve is the locally irregular
sequence
\[
f_t=\log\{1+\exp(z_t)\},\qquad t\in\mathcal I_d,
\]
where $(z_t)$ is obtained by smoothing Gaussian white noise with a
resolution-dependent moving-average filter. This construction produces
heterogeneous local regularity without imposing a global Sobolev-type
structure. The target curve is generated once for each grid size and held
fixed across replications. The main results are reported for
$d=100$.

We generate observations according to
\[
Y=f+0.055\,\varepsilon,
\qquad
\operatorname{Cov}(\varepsilon_j,\varepsilon_k)
=\rho_0^{|j-k|},
\qquad \rho_0=0.6.
\]
The innovations used to generate $\varepsilon$ are standardized Gaussian,
Laplace, or Student-$t$ variables with 3 degrees of freedom. Thus, the three settings share the same
AR(1) covariance structure but differ in their marginal tail behaviour.

We compare the covariance-adapted sequential smoother, using orders
$r=1,2,3$, under both oracle and plug-in covariance geometries; this is the
sequential covariance-adapted analogue of the construction in
Section~\ref{sec:3}, whereas the plug-in theory of Section~\ref{sec:5} is
stated for the simultaneous smoother $S_\vartheta$ of
\eqref{eq:cov-adapted-smoother}. For the plug-in version, the AR(1)
parameter is estimated from an independent
calibration sample of size 200 generated using the same innovation family
as the corresponding reconstruction experiment. To separate covariance
adaptation from the underlying uncorrelated difference construction, we
additionally include the same smoother constructed under the deliberately
misspecified white-noise geometry $\rho=0$ (the independence-based
construction of the first experiment), as well as a conventional sequential
smoother using ordinary differences of orders $r=1,2,3$.
We further include the two discrete constructions introduced in
Section~\ref{sec:3}. The sequential discrete smoother applies the blended
difference penalties successively over orders $r=1,\ldots,4$, with the
blend parameter fixed at $\eta=0.5$. The convex discrete smoother uses the
second-order blended penalty, with $\eta$ jointly optimized to determine
the balance between conventional and uncorrelated differences. Finally,
we include Fourier, B-spline, and Gaussian kernel smoothers as conventional
basis and kernel competitors. For the basis methods, the number of basis 
functions is set sufficiently large so that smoothness is 
primarily controlled through penalization.

To make the comparison concern reconstruction rather than tuning
performance, the regularization parameters are selected by oracle
minimization of the reconstruction MSE against the known target. For the
sequential methods, a common multistart optimization scheme is used to
reduce sensitivity to local minima. Results are averaged over 100 Monte
Carlo replications and reported in \autoref{tab:reconstruction}.

\begin{table}[t]
\caption{Reconstruction MSE and SD (both multiplied by $100$) for the
locally irregular curve with $d=100$ under AR(1) covariance with
$\rho_0=0.6$. The covariance parameter of the plug-in method is estimated
from an independent calibration sample of size 200. Results are averaged
over 100 replications.}
\label{tab:reconstruction}
\centering
\begin{tabular}{@{}lccc@{}}
\hline
Method & Gaussian & Laplace & Student-$t_3$\\
\hline
Covariance-adapted (oracle)
& 0.239 (0.052)
& 0.229 (0.063)
& 0.207 (0.085)\\
Covariance-adapted (plug-in)
& 0.239 (0.052)
& 0.229 (0.063)
& 0.207 (0.085)\\
Uncorrelated geometry (white)
& 0.240 (0.052)
& 0.230 (0.063)
& 0.207 (0.085)\\
Conventional sequential ($R=3$)
& 0.250 (0.053)
& 0.239 (0.066)
& 0.215 (0.087)\\
Sequential discrete
& 0.243 (0.052)
& 0.234 (0.067)
& 0.209 (0.087)\\
Convex discrete
& 0.266 (0.055)
& 0.256 (0.069)
& 0.231 (0.087)\\
Fourier
& 0.291 (0.053)
& 0.280 (0.070)
& 0.255 (0.084)\\
B-spline
& 0.301 (0.054)
& 0.291 (0.070)
& 0.266 (0.084)\\
Gaussian kernel
& 0.264 (0.054)
& 0.253 (0.071)
& 0.226 (0.091)\\
\hline
\end{tabular}
\end{table}

The oracle and plug-in covariance-adapted smoothers yield virtually
identical reconstruction errors under all three innovation distributions,
showing that estimation of the covariance parameter introduces negligible
additional reconstruction error in this experiment. Both also improve on
the conventional sequential $R=3$ smoother and on the conventional basis
and kernel competitors.

The construction based on the white-noise geometry performs very similarly
to the covariance-adapted versions. Thus, in this reconstruction setting,
most of the improvement over the conventional sequential construction is
attributable to the decorrelation and normalization of the difference
operators, whereas adaptation to the specific AR(1) covariance geometry
provides a smaller additional gain. This is consistent with the still-small
absolute magnitude of the operator-level misspecification discrepancy at
$\rho=0.6$ in \autoref{tab:geometry}, which becomes substantially larger
at higher values of $\rho$. The same qualitative pattern persists under
Laplace and Student-$t$ innovations, indicating that the qualitative
reconstruction ranking is stable across these departures from Gaussian
marginals. The Gaussian setting corresponds directly to the efficiency
theory of Section~\ref{sec:5}; the Laplace and Student-$t$ settings use the
same moment conditions but are included only as robustness experiments
outside that theory, without asserting an unproved consistency or
efficiency result for non-Gaussian innovations.

\appendix
\titleformat{\section}[hang]{\large\normalfont\scshape\filcenter}{Appendix~\thesection.}{1em}{}

\section{Proofs of formal statements}\label{appn}

\begin{proof}[Proof of Proposition~\ref{prop:1}]
For $q\geq 0$, let
\[
p_q(s)=s^q.
\]
Applying the stencil to $p_q$ at the origin gives
\[
(\bar D^{(r)}p_q)(0)
=
\sum_\ell w_\ell^{(r)}p_q(\ell)
=
\sum_\ell \ell^q w_\ell^{(r)}
=
M_q(w^{(r)}).
\]
Thus the discrete moments can be obtained by applying the operator to
monomials and evaluating the result at the origin.
We first consider $q<r$. Each application of
the first-difference operator $\Delta^{(1)}=1-T^{-1}$ reduces the degree of a nonconstant polynomial by one.
Consequently,
\[
\Delta^{(r)}p_q=0,
\qquad q=0,\ldots,r-1.
\]
Since $\bar D^{(r)}=W_r(T)\Delta^{(r)}$, it follows that
\[
\bar D^{(r)}p_q
=
W_r(T)\Delta^{(r)}p_q
=
0.
\]
Evaluating at the origin therefore gives
\[
M_q(w^{(r)})
=
(\bar D^{(r)}p_q)(0)
=
0,
\qquad q=0,\ldots,r-1.
\]
For $q=r$, the $r$th conventional difference of $p_r(s)=s^r$ is the
constant
\[
\Delta^{(r)}p_r=r!.
\]
Hence
\[
\bar D^{(r)}p_r
=
W_r(T)(r!)
=
\sum_{j\in J_r}a_j^{(r)}T^j(r!).
\]
Translations leave a constant sequence unchanged, so
$T^j(r!)=r!$ for every $j$. Therefore,
\[
M_r(w^{(r)})
=
(\bar D^{(r)}p_r)(0)
=
r!\sum_{j\in J_r}a_j^{(r)}.
\]
Thus all moments below order $r$ vanish, while the moment of order $r$
is nonzero whenever $\sum_{j\in J_r}a_j^{(r)}\neq0$, proving that the
stencil has exact difference order $r$.
We next establish the derivative approximation. Let
\[
t_i=t_0+ih,
\qquad
x_i=f(t_i),
\]
so that, because the grid is equally spaced,
\[
t_{i+\ell}=t_i+\ell h.
\]
Applying the stencil to the sampled function therefore gives
\[
(\bar D^{(r)}x)_i
=
\sum_\ell w_\ell^{(r)}f(t_i+\ell h).
\]
For each fixed stencil offset $\ell$, Taylor expansion about $t_i$ gives,
for sufficiently small $h$,
\[
f(t_i+\ell h)
=
\sum_{q=0}^{r}
\frac{(\ell h)^q}{q!}f^{(q)}(t_i)
+
O(h^{r+1}).
\]
Since the stencil has fixed finite support, multiplying by
$w_\ell^{(r)}$ and summing over $\ell$ yields
\[
\begin{aligned}
\sum_\ell w_\ell^{(r)}f(t_i+\ell h)
&=
\sum_{q=0}^{r}
\frac{h^q}{q!}f^{(q)}(t_i)
\sum_\ell \ell^q w_\ell^{(r)}
+
O(h^{r+1})\\
&=
\sum_{q=0}^{r}
\frac{h^q}{q!}f^{(q)}(t_i)
M_q(w^{(r)})
+
O(h^{r+1}).
\end{aligned}
\]
The moment identities already proved eliminate all terms with $q<r$.
Hence
\[
\sum_\ell w_\ell^{(r)}f(t_i+\ell h)
=
\frac{h^r}{r!}
M_r(w^{(r)})f^{(r)}(t_i)
+
O(h^{r+1}).
\]
Substituting the expression for $M_r(w^{(r)})$ gives
\[
\sum_\ell w_\ell^{(r)}f(t_i+\ell h)
=
h^r
\left(\sum_{j\in J_r}a_j^{(r)}\right)
f^{(r)}(t_i)
+
O(h^{r+1}),
\]
which proves the stated approximation.
Finally, suppose that
\[
w_{-\ell}^{(r)}=(-1)^r w_\ell^{(r)}.
\]
For every $\ell>0$, the contributions at $\ell$ and $-\ell$ to the
$(r+1)$st moment cancel:
\[
\begin{aligned}
w_\ell^{(r)}\ell^{r+1}
+
w_{-\ell}^{(r)}(-\ell)^{r+1}
&=
w_\ell^{(r)}\ell^{r+1}
+
(-1)^r w_\ell^{(r)}(-1)^{r+1}\ell^{r+1}\\
&=0.
\end{aligned}
\]
The term at $\ell=0$ also contributes zero. Therefore,
\[
M_{r+1}(w^{(r)})=0.
\]
Expanding Taylor's formula one order further, the term of order
$h^{r+1}$ is
\[
\frac{h^{r+1}}{(r+1)!}
f^{(r+1)}(t_i)M_{r+1}(w^{(r)}),
\]
which vanishes. Thus, for $f\in C^{r+2}$, the remainder improves to
$O(h^{r+2})$.
\end{proof}

\begin{proof}[Proof of Corollary~\ref{col:1}]
After zero-padding stencils of different widths to a common support,
independence and unit variance of the entries of $\varepsilon$ imply that the covariance
between two transformed variables at the same location equals
$\sum_\ell w_\ell^{(r)}w_\ell^{(s)}$. This is zero by the decorrelation
construction and equals one when $r=s$ by normalization. The expectation is
zero by linearity. If $\varepsilon$ is Gaussian, the transformed variables are jointly
Gaussian linear functionals of $\varepsilon$, so zero covariance implies independence.
\end{proof}

\begin{proof}[Proof of Proposition~\ref{prop:spectral-normalization}]
Unit-variance normalization gives
$\sum_\ell |w_\ell^{(r)}|^2=1$, and the conclusion follows immediately from
Parseval's identity.
\end{proof}

\begin{proof}[Proof of Proposition~\ref{prop:multiorder}]
The result follows from $Z(t)\sim\mathcal N(0,I_R)$, since its squared
components are independent $\chi_1^2$ random variables.
\end{proof}

\begin{proof}[Proof of Proposition~\ref{prop:variance-optimality}]
Expanding the quadratic-form variance gives
\[
\operatorname{tr}\{(AG)^2\}
=
\sum_{r,s=1}^R\alpha_r\alpha_sG_{rs}^2
=
\sum_{r=1}^R\alpha_r^2
+
\sum_{r\neq s}\alpha_r\alpha_sG_{rs}^2,
\]
where $G_{rr}=1$. Since $\alpha_r>0$, the second term is non-negative and
vanishes if and only if all off-diagonal entries of $G$ vanish.
\end{proof}

\begin{proof}[Proof of Proposition~\ref{prop:cov-consistency}]
By construction,
$w_\vartheta^{(r)}\in\mathcal V_r^\pm(J_r)\subseteq\mathcal V_r(J_r)$ and
therefore remains a finite linear combination of translates of
$\Delta^{(r)}$. Proposition~\ref{prop:1} relies only on the algebraic
identities
$\Delta^{(r)}p_q\equiv0$ for $q<r$ and
$\Delta^{(r)}p_r\equiv r!$, together with the parity of the resulting
stencil. None of these properties depends on the reference covariance
$\Sigma_\vartheta$.
\end{proof}

\begin{proof}[Proof of Proposition~\ref{prop:cov-smoothness}]
The argument proceeds by induction on $r$. The initial vector
$w_\vartheta^{(0)}=e_0$ is constant. Suppose that
$w_\vartheta^{(0)},\ldots,w_\vartheta^{(r-1)}$ are continuously
differentiable on $U$. Since $\Sigma_\vartheta$ is continuously
differentiable, so is $C_r(\vartheta)$.

The condition $\dim\ker C_r(\vartheta)=1$ implies
$\operatorname{rank}C_r(\vartheta)=m_r-1$. Hence, locally around any point
of $U$, a nonzero $(m_r-1)\times(m_r-1)$ minor can be selected. Solving the
corresponding homogeneous system after fixing one free coordinate gives a
nonzero kernel representative $v_r(\vartheta)$ whose entries are rational
functions of the entries of $C_r(\vartheta)$ and are therefore continuously
differentiable locally.

Since $\Sigma_\vartheta\succ0$,
$v_r(\vartheta)^\top B_r^\top\Sigma_\vartheta B_rv_r(\vartheta)>0$.
Thus
\[
a_r(\vartheta)
=
\frac{v_r(\vartheta)}
{\left\{
v_r(\vartheta)^\top
B_r^\top\Sigma_\vartheta B_r
v_r(\vartheta)
\right\}^{1/2}}
\]
is continuously differentiable. The prescribed sign convention is locally
constant wherever its sign-fixing coordinate is nonzero, completing the
induction.
\end{proof}

\begin{proof}[Proof of Lemma~\ref{lem:crossparity}]
Since $J^\top=J$ and $J^2=I$,
\[
u^\top\Sigma_\vartheta v
=
(Ju)^\top\Sigma_\vartheta(Jv)
=
u^\top\Sigma_\vartheta(-v)
=
-u^\top\Sigma_\vartheta v.
\]
Hence $u^\top\Sigma_\vartheta v=0$.
\end{proof}

\begin{proof}[Proof of Corollary~\ref{cor:cov-geometry}]
Let $W_{\vartheta,t}$ denote the matrix whose columns are the zero-padded
stencils $w_\vartheta^{(1)},\ldots,w_\vartheta^{(R)}$ shifted to location
$t$, that is, $W_{\vartheta,t}=(P_tw_\vartheta^{(1)},\ldots,P_tw_\vartheta^{(R)})$.
Then $Z_\vartheta(t)=W_{\vartheta,t}^\top\varepsilon$, and by translation
invariance
$(P_tw_\vartheta^{(r)})^\top\Sigma_\vartheta(P_tw_\vartheta^{(s)})
=w_\vartheta^{(r)\top}\Sigma_\vartheta w_\vartheta^{(s)}$
for every interior $t$, the same value used to construct $w_\vartheta^{(r)}$.
Hence
\[
\operatorname{cov}_{Q_\vartheta}\{Z_\vartheta(t)\}
=
W_{\vartheta,t}^\top\Sigma_\vartheta W_{\vartheta,t}
=
I_R
\]
by the normalization and mutual orthogonality imposed in the construction.
The conclusion follows from Gaussianity.
\end{proof}

\begin{proof}[Proof of Proposition~\ref{prop:S-smooth}]
Each $\bar D_\vartheta^{(r)}$ is obtained by fixed linear reindexing of
$w_\vartheta^{(r)}$ and is therefore $C^1$ on $U$ by
Proposition~\ref{prop:cov-smoothness}. Consequently $\mathfrak P_\vartheta$
is $C^1$ on $U$. Since $\mathfrak P_\vartheta\succeq0$,
$I_d+\mathfrak P_\vartheta$ is invertible, and matrix inversion is
continuously differentiable on the set of invertible matrices.
Differentiating $S_\vartheta(I_d+\mathfrak P_\vartheta)=I_d$ on $U$ gives the
first identity in \eqref{eq:Sderivative}; the second follows from the
product rule.
\end{proof}

\begin{proof}[Proof of Proposition~\ref{prop:Dtau}]
Write $w=w_\tau^{(r)}$, so that $w^\top\Sigma_\tau w=1$ and
$h_{\tau,r}(x)=(w^\top x)^2-1$. The log-density of $\mathcal N(0,\Sigma_\vartheta)$
at a fixed $x$ is $-\tfrac12\log|\Sigma_\vartheta|-\tfrac12x^\top\Sigma_\vartheta^{-1}x$
up to an additive constant, so, using
\[
\left.\partial_{\vartheta_j}\log|\Sigma_\vartheta|\right|_{\vartheta=\tau}
=
\operatorname{tr}(B_{\tau,j}\Sigma_\tau)
\]
and
\[
\left.\partial_{\vartheta_j}\Sigma_\vartheta^{-1}\right|_{\vartheta=\tau}
=
-B_{\tau,j},
\]
where
\[
B_{\tau,j}:=\Sigma_\tau^{-1}\dot\Sigma_{\tau,j}\Sigma_\tau^{-1},
\]
the score has the
explicit form
\[
\kappa_{\tau,j}(x)
=
\tfrac12\left\{x^\top B_{\tau,j}x-\operatorname{tr}(B_{\tau,j}\Sigma_\tau)\right\}.
\]
By the Gaussian fourth-moment identity, for $\varepsilon\sim\mathcal N(0,\Sigma_\tau)$
and symmetric matrices $A_1,A_2$,
\[
E_\tau\{(\varepsilon^\top A_1\varepsilon)(\varepsilon^\top A_2\varepsilon)\}
=
\operatorname{tr}(A_1\Sigma_\tau)\operatorname{tr}(A_2\Sigma_\tau)
+2\operatorname{tr}(A_1\Sigma_\tau A_2\Sigma_\tau).
\]
Applying this with $A_1=ww^\top$, $A_2=B_{\tau,j}$, and using
$\Sigma_\tau B_{\tau,j}\Sigma_\tau=\dot\Sigma_{\tau,j}$ together with
$w^\top\Sigma_\tau w=1$, gives
\[
E_\tau\{(w^\top\varepsilon)^2(\varepsilon^\top B_{\tau,j}\varepsilon)\}
=
\operatorname{tr}(B_{\tau,j}\Sigma_\tau)+2w^\top\dot\Sigma_{\tau,j}w.
\]
Combining this with $E_\tau\{(w^\top\varepsilon)^2\}=1$ and
$E_\tau\{\varepsilon^\top B_{\tau,j}\varepsilon\}=\operatorname{tr}(B_{\tau,j}\Sigma_\tau)$
yields
\[
E_\tau\{h_{\tau,r}(\varepsilon)\kappa_{\tau,j}(\varepsilon)\}
=
\tfrac12
E_\tau\big[
\{(w^\top\varepsilon)^2-1\}
\{\varepsilon^\top B_{\tau,j}\varepsilon-\operatorname{tr}(B_{\tau,j}\Sigma_\tau)\}
\big]
=
w^\top\dot\Sigma_{\tau,j}w,
\]
which gives \eqref{eq:Dtau}. Equivalently, differentiating
$w_\vartheta^{(r)\top}\Sigma_\vartheta w_\vartheta^{(r)}\equiv1$
shows directly that
\[
E_\tau\!\left\{
\left.\partial_j h_{\vartheta,r}(\varepsilon)\right|_{\vartheta=\tau}
\right\}
=
2\dot w_{\tau,j}^{(r)\top}\Sigma_\tau w_\tau^{(r)}
=
-w_\tau^{(r)\top}\dot\Sigma_{\tau,j}w_\tau^{(r)}
=
\Gamma_{\tau,rj}.
\]
\end{proof}

\begin{proof}[Proof of Proposition~\ref{prop:Omega}]
Write
$h_{\tau,r}(\varepsilon)=\varepsilon^\top A_r\varepsilon-1$, where
$A_r=w_\tau^{(r)}w_\tau^{(r)\top}$. For a centered Gaussian vector,
\[
\operatorname{cov}
(\varepsilon^\top A_r\varepsilon,\varepsilon^\top A_s\varepsilon)
=
2\operatorname{tr}
(A_r\Sigma_\tau A_s\Sigma_\tau).
\]
Since $A_r$ and $A_s$ have rank one, the right-hand side equals
$2(w_\tau^{(r)\top}\Sigma_\tau w_\tau^{(s)})^2$,
which proves \eqref{eq:Omega}.
\end{proof}

\begin{proof}[Proof sketch for Theorem~\ref{thm:plugin}]
By Proposition~\ref{prop:cov-smoothness},
$\vartheta\mapsto w_\vartheta^{(r)}$ is locally continuous for every
$r$. The conditions for Hellinger continuity in Lemma~2.5 of
\citet{Schick:2001} are therefore verified as follows. First, for each
fixed $x$, $h_\vartheta(x)s_\vartheta(x)$ is continuous in $\vartheta$,
where $s_\vartheta=f_\vartheta^{1/2}$, since both
$\Sigma_\vartheta$ and $w_\vartheta^{(1)},\ldots,w_\vartheta^{(R)}$
vary continuously. Hence
$h_\vartheta s_\vartheta\to h_{\vartheta_0}s_{\vartheta_0}$ in measure
on sets of finite measure. Second, the moments
\[
\int h_\vartheta h_\vartheta^\top f_\vartheta\,d\mu
\]
are finite-dimensional Gaussian moments of degree at most four and are
continuous in $\vartheta$ under the stated assumptions. Thus
$\vartheta\mapsto h_\vartheta$ is Hellinger continuous at
$\vartheta_0$. Theorem~2.3 of \citet{Schick:2001} then gives
\eqref{eq:plugin-H} with coefficient matrix
\[
-\int
h_{\vartheta_0}
\kappa_{\vartheta_0}^{\top}
\,dF_{\vartheta_0},
\]
where $\kappa_{\vartheta_0}$ is the Hellinger derivative in the
terminology of \citet{Schick:2001}; for the present smoothly parametrized
Gaussian family, it is the ordinary likelihood score. Hence this matrix is
\[
-E_{\vartheta_0}
\{h_{\vartheta_0}(\varepsilon)
\kappa_{\vartheta_0}(\varepsilon)^\top\}
=
\Gamma_{\vartheta_0},
\]
by Proposition~\ref{prop:Dtau}, and \eqref{eq:plugin-H} follows.
\end{proof}

\begin{proof}[Proof of Corollary~\ref{cor:vartheta-normal}]
Hansen's asymptotic-normality theorem \citep[Thm.~3.1]{Hansen:1982} is
stated for a $k\times R$ limiting weighting matrix $a_0^*$, giving
\[
\sqrt n(\widehat\vartheta_n-\vartheta_0)
\Rightarrow
\mathcal N\{0,
(a_0^*\Gamma_{\vartheta_0})^{-1}a_0^*S_wa_0^{*\top}(a_0^*\Gamma_{\vartheta_0})^{-1\top}
\},
\]
where $S_w$ denotes the long-run variance of $h_{\vartheta_0}(\varepsilon)$
along the calibration sequence. For an estimator constructed by minimizing
the quadratic criterion $H_n(\vartheta)^\top\mathcal W H_n(\vartheta)$,
Hansen's own reduction of the linear-combination form to this weighted
form, via his Eqs.~(6)--(7), gives $a_0^*=\Gamma_{\vartheta_0}^\top\mathcal
W$. Substituting this into the covariance above, and using that
$\Gamma_{\vartheta_0}^\top\mathcal W\Gamma_{\vartheta_0}$ is symmetric,
yields exactly
\[
(\Gamma_{\vartheta_0}^\top\mathcal W\Gamma_{\vartheta_0})^{-1}
\Gamma_{\vartheta_0}^\top\mathcal WS_w\mathcal W\Gamma_{\vartheta_0}
(\Gamma_{\vartheta_0}^\top\mathcal W\Gamma_{\vartheta_0})^{-1}.
\]
Since the calibration sample $\varepsilon_1,\ldots,\varepsilon_n$ is
independent and identically distributed, this is Hansen's Case~(i)
\citep[Sec.~3]{Hansen:1982}, under which the long-run variance collapses
to the ordinary covariance of a single term,
$S_w=\operatorname{cov}_{\vartheta_0}\{h_{\vartheta_0}(\varepsilon)\}=\Omega_{\vartheta_0}$
by Proposition~\ref{prop:Omega}. This gives $V_{\vartheta_0}(\mathcal W)$
exactly as stated, with population Jacobian $\Gamma_{\vartheta_0}$ by
Proposition~\ref{prop:Dtau}.
\end{proof}

\begin{proof}[Proof of Proposition~\ref{prop:fhat-expansion}]
By Proposition~\ref{prop:S-smooth} and the first-order expansion of the
matrix-valued map $\vartheta\mapsto S_\vartheta$,
\[
S_{\widehat\vartheta_n}-S_{\vartheta_0}
=
\sum_{j=1}^k
\dot S_{\vartheta_0,j}
(\widehat\vartheta_{n,j}-\vartheta_{0,j})
+o_p(n^{-1/2})
\]
in matrix norm. Since $\|X\|=O_p(1)$, multiplication by $X$ preserves the
$o_p(n^{-1/2})$ remainder, and \eqref{eq:Sderivative} yields
\eqref{eq:fhat-expansion}.
\end{proof}

\begin{proof}[Proof of Theorem~\ref{thm:efficient-plugin-risk}]
Proposition~\ref{prop:fhat-expansion} gives, for fixed $x$,
\[
\sqrt n
\left\{
\widehat f_{\widehat\vartheta_n}(x)
-
\widehat f_{\vartheta_0}(x)
\right\}
=
L(x)\sqrt n(\widehat\vartheta_n-\vartheta_0)+o_p(1).
\]
By Corollary~\ref{cor:vartheta-normal} and the continuous mapping theorem,
the right-hand side converges in distribution to $L(x)Z$, where
$Z\sim\mathcal N\{0,V_{\vartheta_0}(\mathcal W)\}$, and hence
$n\|\widehat f_{\widehat\vartheta_n}(x)-\widehat f_{\vartheta_0}(x)\|^2
\Rightarrow\|L(x)Z\|^2$. Since this sequence is uniformly integrable by
assumption, convergence in distribution upgrades to convergence of the
corresponding means, and $E\|L(x)Z\|^2=\operatorname{tr}\{L(x)^\top
L(x)V_{\vartheta_0}(\mathcal W)\}$ gives \eqref{eq:conditional-plugin-risk}.
By Corollary~\ref{cor:efficiency}, $\mathcal W=\Omega_{\vartheta_0}^{-1}$ is
the efficient GMM weight for the moment conditions $h_\vartheta$; in the
notation of Corollary~\ref{cor:vartheta-normal}, this is exactly the
Loewner-order statement
\[
V_{\vartheta_0}(\mathcal W)
\succeq
V_{\vartheta_0}(\Omega_{\vartheta_0}^{-1})
\]
for every admissible $\mathcal W$. Since $L(x)^\top L(x)\succeq0$,
trace monotonicity under the Loewner order yields
\[
\operatorname{tr}
\left\{
L(x)^\top L(x)V_{\vartheta_0}(\mathcal W)
\right\}
\ge
\operatorname{tr}
\left\{
L(x)^\top L(x)
V_{\vartheta_0}(\Omega_{\vartheta_0}^{-1})
\right\}.
\]
The efficient weight does not depend on $x$ or
$\boldsymbol\alpha$, so the inequality holds simultaneously for all such
choices. Under decorrelation,
$\Omega_{\vartheta_0}^{-1}=\tfrac12I_R$, which is equivalent in GMM to the
unweighted choice $\mathcal W=I_R$.
\end{proof}

\begin{proof}[Proof of Corollary~\ref{cor:integrated-plugin-risk}]
By Proposition~\ref{prop:fhat-expansion}, independence of $X$ and the
calibration sample, and Corollary~\ref{cor:vartheta-normal},
\[
\sqrt n
\left\{
\widehat f_{\widehat\vartheta_n}(X)
-
\widehat f_{\vartheta_0}(X)
\right\}
\Rightarrow
L(X)Z,
\]
where $Z\sim\mathcal N\{0,V_{\vartheta_0}(\mathcal W)\}$ is independent of
$X$. The assumed uniform integrability therefore gives
\[
nE
\left\|
\widehat f_{\widehat\vartheta_n}(X)
-
\widehat f_{\vartheta_0}(X)
\right\|^2
\longrightarrow
E\|L(X)Z\|^2.
\]
By independence,
\[
E\|L(X)Z\|^2
=
\operatorname{tr}
\left\{
E[L(X)^\top L(X)]
V_{\vartheta_0}(\mathcal W)
\right\},
\]
which gives \eqref{eq:integrated-plugin-risk}. The optimality statement
follows from the same Loewner-order argument as in
Theorem~\ref{thm:efficient-plugin-risk}.
\end{proof}

\section{Additional simulations}\label{appB}

\paragraph{\textbf{Smooth sinusoidal curve.}}
As a complementary reconstruction experiment, we replace the locally
irregular target of Section~\ref{sec:7} by the globally smooth sinusoid
\[
f_t=\sin(6\pi x_t),\qquad t\in\mathcal I_d,
\]
where $(x_t)_{t\in\mathcal I_d}$ is an equally spaced grid on $[0,1]$, and
consider resolutions $d\in\{50,100,200\}$. All other aspects of the
reconstruction experiment are unchanged: observations are generated under
AR(1) covariance with $\rho_0=0.6$ and noise scale $0.055$, using standardized
Gaussian, Laplace, or Student-$t_3$ innovations.

Since the effect of covariance adaptation and its estimation is examined
directly in Section~\ref{sec:7}, we omit those variants here and focus on the
comparison between the discrete smoothing constructions and the conventional
competitors. We compare the conventional sequential smoother with $R=3$,
the sequential discrete smoother with $R=4$ and $\eta=0.5$, the convex
discrete smoother, and the Fourier, B-spline, and Gaussian kernel smoothers.
Regularization parameters are selected by oracle minimization of reconstruction
MSE, using the same multistart protocol for the sequential methods. Results
are averaged over 100 Monte Carlo replications.

The smooth-signal experiment complements the locally irregular reconstruction
experiment of Section~\ref{sec:7}. The sequential discrete smoother performs
particularly well, attaining the lowest or essentially lowest reconstruction
error across resolutions and innovation distributions and improving on the
Fourier, B-spline, and Gaussian kernel competitors. Its advantage becomes
clearer as the grid resolution increases, showing that the discrete
construction remains competitive also for a globally smooth target.

\begin{table}[t!]
\caption{Reconstruction MSE and SD (both multiplied by $100$) for the
smooth sinusoidal curve under AR(1) covariance with $\rho_0=0.6$ at different
grid resolutions. Results are averaged over 100 replications.}
\label{tab:smoothing_sinusoid}
\centering
\begin{tabular}{@{}lccc@{}}
\hline
Method & Gaussian & Laplace & Student-$t_3$\\
\hline

\multicolumn{4}{l}{\emph{$d=50$}}\\[2pt]
Conventional sequential ($R=3$)
& \textbf{0.195 (0.068)}
& \textbf{0.218 (0.100)}
& \textbf{0.162 (0.094)}\\
Sequential discrete
& \textbf{0.195 (0.069)}
& \textbf{0.218 (0.100)}
& \textbf{0.162 (0.094)}\\
Convex discrete
& 0.228 (0.072)
& 0.253 (0.102)
& 0.204 (0.098)\\
Fourier
& 0.201 (0.070)
& 0.223 (0.103)
& 0.170 (0.100)\\
B-spline
& 0.216 (0.070)
& 0.238 (0.104)
& 0.180 (0.101)\\
Gaussian kernel
& 0.271 (0.074)
& 0.294 (0.111)
& 0.228 (0.119)\\[5pt]

\multicolumn{4}{l}{\emph{$d=100$}}\\[2pt]
Conventional sequential ($R=3$)
& 0.133 (0.050)
& 0.122 (0.057)
& 0.112 (0.083)\\
Sequential discrete
& \textbf{0.127 (0.049)}
& \textbf{0.117 (0.056)}
& \textbf{0.107 (0.080)}\\
Convex discrete
& 0.164 (0.063)
& 0.162 (0.074)
& 0.170 (0.093)\\
Fourier
& 0.145 (0.049)
& 0.134 (0.064)
& 0.125 (0.089)\\
B-spline
& 0.156 (0.053)
& 0.145 (0.064)
& 0.134 (0.094)\\
Gaussian kernel
& 0.226 (0.058)
& 0.215 (0.076)
& 0.196 (0.118)\\[5pt]

\multicolumn{4}{l}{\emph{$d=200$}}\\[2pt]
Conventional sequential ($R=3$)
& 0.075 (0.029)
& 0.074 (0.030)
& 0.066 (0.037)\\
Sequential discrete
& \textbf{0.074 (0.028)}
& \textbf{0.072 (0.030)}
& \textbf{0.064 (0.036)}\\
Convex discrete
& 0.097 (0.031)
& 0.097 (0.035)
& 0.129 (0.065)\\
Fourier
& 0.088 (0.031)
& 0.084 (0.030)
& 0.076 (0.042)\\
B-spline
& 0.097 (0.031)
& 0.094 (0.032)
& 0.083 (0.044)\\
Gaussian kernel
& 0.166 (0.038)
& 0.163 (0.039)
& 0.146 (0.066)\\
\hline
\end{tabular}
\end{table}

\newpage
\bibliographystyle{apalike-dashed}
 {\footnotesize
\bibliography{ref}}

@incollection{Mizuta:2023,
  author    = {Mizuta, Masahiro},
  title     = {Discrete Functional Data Analysis Based on Discrete Difference},
  booktitle = {Analysis of Categorical Data from Historical Perspectives: Essays in Honour of Shizuhiko Nishisato},
  editor    = {Beh, Eric J. and Lombardo, Rosaria and Clavel, Jose G.},
  publisher = {Springer Nature Singapore},
  address   = {Singapore},
  year      = {2023},
  pages     = {487--492},
  doi       = {10.1007/978-981-99-5329-5_28},
}

@incollection{Mizuta:2006,
  author    = {Mizuta, Masahiro},
  title     = {Discrete Functional Data Analysis},
  booktitle = {Compstat 2006 -- Proceedings in Computational Statistics},
  editor    = {Rizzi, Alfredo and Vichi, Maurizio},
  publisher = {Physica-Verlag HD},
  address   = {Heidelberg},
  year      = {2006},
  pages     = {361--369},
  isbn      = {978-3-7908-1709-6},
  doi       = {10.1007/978-3-7908-1709-6_28}
}

@article{Hall:Opsomer:2005,
  author = {Hall, Peter and Opsomer, J. D.},
  title = {Theory for penalised spline regression},
  journal = {Biometrika},
  volume = {92},
  number = {1},
  pages = {105--118},
  year = {2005},
  doi={10.1093/biomet/92.1.105}
 }

@article{Mueller:1987,
  author  = {M{\"u}ller, Hans-Georg},
  title   = {Weighted local regression and kernel methods for nonparametric curve fitting},
  journal = {J. Am. Stat. Assoc.},
  volume  = {82},
  number  = {397},
  pages   = {231--238},
  year    = {1987},
  doi = {10.1080/01621459.1987.10478425}
}

@article{Schick:2001,
  author  = {Schick, Anton},
  title   = {On asymptotic differentiability of averages},
  journal = {Statist. Probab. Lett.},
  volume  = {51},
  number  = {1},
  pages   = {15--23},
  year    = {2001},
  issn    = {0167-7152},
  doi     = {10.1016/S0167-7152(00)00132-2},
}

@article{Belhakem:Picard:Rivoirard:Roche:2025,
  author  = {Belhakem, Rachid and Picard, Fran{\c{c}}ois and Rivoirard, Vincent and Roche, Alexis},
  title   = {Minimax Estimation of Functional Principal Components from Noisy Discretized Functional Data},
  journal = {Scand. J. Stat.},
  volume  = {52},
  number  = {1},
  pages   = {38--80},
  year    = {2025},
  doi     = {10.1111/sjos.12719},
}

@article{Whittaker:1923,
  author  = {Whittaker, E. T.},
  title   = {On a New Method of Graduation},
  journal = {Proc. Edinb. Math. Soc.},
  volume  = {41},
  pages   = {63--75},
  year    = {1922},
  doi     = {10.1017/S0013091500077853},
}

@article{Hodrick:Prescott:1997,
  author  = {Hodrick, Robert J. and Prescott, Edward C.},
  title   = {Postwar {U.S.} Business Cycles: An Empirical Investigation},
  journal = {J. Money Credit Bank.},
  volume  = {29},
  number  = {1},
  pages   = {1--16},
  year    = {1997},
  doi     = {10.2307/2953682},
}

@article{Eilers:2003,
  author  = {Eilers, Paul H. C.},
  title   = {A Perfect Smoother},
  journal = {Anal. Chem.},
  volume  = {75},
  number  = {14},
  pages   = {3631--3636},
  year    = {2003},
  doi     = {10.1021/ac034173t},
}

@article{Claeskens:Krivobokova:Opsomer:2009,
  author  = {Claeskens, Gerda and Krivobokova, Tatyana and Opsomer, Jean D.},
  title   = {Asymptotic Properties of Penalized Spline Estimators},
  journal = {Biometrika},
  volume  = {96},
  number  = {3},
  pages   = {529--544},
  year    = {2009},
  doi     = {10.1093/biomet/asp035},
}

@article{Tibshirani:2022,
  author  = {Tibshirani, Ryan J.},
  title   = {Divided Differences, Falling Factorials, and Discrete Splines:
             Another Look at Trend Filtering and Related Problems},
  journal = {Foundations and Trends in Machine Learning},
  volume  = {15},
  number  = {6},
  pages   = {694--846},
  year    = {2022},
  doi     = {10.1561/2200000099},
}

@article{Eilers:Marx:Durban:2015,
  author = {Eilers, Paul H. C. and Marx, Brian D. and Durb{\'a}n, Mar{\'i}a},
  title = {Twenty years of {P}-splines},
  journal = {SORT. Statistics and Operations Research Transactions},
  volume = {39},
  number = {2},
  pages = {149--186},
  year = {2015},
  url = {https://raco.cat/index.php/SORT/article/view/302258.}
}

@Book{vanderVaart:1998,
  Title     = {Asymptotic Statistics},
  Author    = {van der Vaart, A. W.},
  Publisher = {Cambridge University Press},
  Address   = {Cambridge},
  Year      = {1998},
}

@article{Hansen:1982,
  author  = {Hansen, Lars Peter},
  title   = {Large Sample Properties of Generalized Method of Moments Estimators},
  journal = {Econometrica},
  volume  = {50},
  number  = {4},
  pages   = {1029--1054},
  year    = {1982},
  doi     = {10.2307/1912775},
}

@article{Vidal:2026,
  author  = {Vidal, Marc},
  title   = {Beyond the Orthogonal Group: A Discussion on Robust {ICA}
             and Sphering Transformations},
  journal = {Advances in Data Analysis and Classification},
  year    = {2026},
  doi     = {10.1007/s11634-026-00699-0},
}

\end{document}